\documentclass[12pt]{amsart}
\usepackage{amssymb}
\usepackage{verbatim}
\usepackage{graphicx}
\usepackage[cmtip,all]{xy}
\usepackage{amsmath}
\usepackage{amsfonts,amssymb}
\usepackage{times}
\usepackage{amscd}
\usepackage{epsfig}
\usepackage{float}
\newtheorem{theo}{Theorem}[section]
\newtheorem{thm}[theo]{Theorem}
\newtheorem{lem}[theo]{Lemma}
\newtheorem{cor}[theo]{Corollary}

\newtheorem{prop}[theo]{Proposition}

\newtheorem*{thmM}{Main Theorem}

\theoremstyle{definition}
\newtheorem{dfn}[theo]{Definition}

\theoremstyle{remark}

\numberwithin{equation}{section}

\def\R{\mathbb{R}}
\def\Z{\mathbb{Z}}
\def\Mc{\mathcal{M}}
\newcommand{\gl}{geolamination{}}
\newcommand{\ql}{quadrilateral{}}

\newcommand{\md}{\mathcal{MD}}
\newcommand{\fup}{\mathcal{P}}
\newcommand{\qcp}{\mathrm{QCP}}

\newcommand{\cm}{\mathrm{CM}}

\newcommand{\coc}{\mathrm{co}}

\newcommand{\cml}{\mathrm{CML}}

\newcommand{\ta}{\mathrm{T\!ag}}

\newcommand{\cp}{\mathcal C}

\newcommand{\cpd}{\mathcal{CPD}}
\newcommand{\cmd}{\mathcal{CMD}}

\newcommand{\rc}{\mathcal{R}}
\newcommand{\zc}{\mathcal{Z}}
\newcommand{\nc}{\mathcal{N}}

\newcommand{\nin}{\notin}

\newcommand{\U}{\mathcal{U}}

\newcommand{\lp}{\mathcal{LP}}

\newcommand{\C}{\mathbb{C}}

\newcommand{\disk}{\mathbb{D}}
\newcommand{\cdisk}{\ol{\mathbb{D}}}

\newcommand{\ol}{\overline}
\newcommand{\ovc}{\ol{c}}

\newcommand{\sm}{\setminus}

\newcommand{\B}{\mathcal{B}}
\newcommand{\Tc}{\mathcal{T}}

\newcommand{\bt}{\ol{t}}
\newcommand{\m}{\ol{m}}
\newcommand{\n}{\ol{n}}
\newcommand{\hell}{\hat{\ell}}

\newcommand{\qml}{\mathrm{QML}}
\newcommand{\bd}{\mathrm{Bd}}

\newcommand{\lam}{\mathcal{L}}
\newcommand{\hlam}{\mathcal{\widehat L}}

\newcommand{\lamm}{\mathcal{L}^m}

\newcommand{\fqcp}{\mathcal{QCP}}

\newcommand{\ch}{\mathrm{CH}}

\newcommand{\si}{\sigma}
\newcommand{\ph}{\varphi}
\newcommand{\uc}{\mathbb{S}}

\newcommand{\om}{\omega}

\newcommand{\g}{\mathfrak{g}}

\newcommand{\mD}{\mathcal{D}}
\newcommand{\Cc}{\mathcal{C}}
\newcommand{\ovk}{\ol{k}}

\def\L{\mathbb{L}}
\def\dc{\mathcal{D}}

\newcommand{\ld}{\L\mathcal{D}}

\renewcommand\le{\leqslant}
\renewcommand\ge{\geqslant}
\def\0{\varnothing}

\begin{document}

\date{January 11, 2015}

\title[Combinatorial models]{Combinatorial models for spaces of\\ cubic polynomials}

\author[A.~Blokh]{Alexander~Blokh}

\thanks{The first and the third named authors were partially
supported by NSF grant DMS--1201450}

\author[L.~Oversteegen]{Lex Oversteegen}

\author[R.~Ptacek]{Ross~Ptacek}

\author[V.~Timorin]{Vladlen~Timorin}

\thanks{The fourth named author was partially supported by
RFBR grants 13-01-12449, 13-01-00969, and AG Laboratory
NRU-HSE, MESRF grant ag. 11.G34.31.0023.}

\address[Alexander~Blokh and Lex~Oversteegen]
{Department of Mathematics\\ University of Alabama at Birmingham\\
Birmingham, AL 35294}

\address[Ross~Ptacek and Vladlen~Timorin]
{Faculty of Mathematics\\
Laboratory of Algebraic Geometry and its Applications\\
Higher School of Economics\\
Vavilova St. 7, 112312 Moscow, Russia }

\address[Vladlen~Timorin]
{Independent University of Moscow\\
Bolshoy Vla-syevskiy Pereulok 11, 119002 Moscow, Russia}

\email[Alexander~Blokh]{ablokh@math.uab.edu}
\email[Lex~Oversteegen]{overstee@math.uab.edu}
\email[Ross~Ptacek]{rptacek@uab.edu}
\email[Vladlen~Timorin]{vtimorin@hse.ru}



\begin{abstract}
A model for the Mandelbrot set is due to Thurston
and is stated in the language of geodesic laminations.
The conjecture that the Mandelbrot set is actually homeomorphic to this model
is equivalent to the celebrated MLC conjecture stating that the Mandelbrot
set is locally connected.
For parameter spaces of higher degree polynomials, even conjectural models
are missing, one possible reason being that the higher degree analog of the
MLC conjecture is known to be false.
We provide a combinatorial model for an essential part of the parameter
space of complex cubic polynomials, namely, for the space of all cubic
polynomials with connected Julia sets all of whose cycles are repelling
(we call such polynomials \emph{dendritic}).
The description of the model turns out to be very similar to that of Thurston.
\end{abstract}

\maketitle

\section{Introduction}\label{s:intro}
The \emph{parameter space} of complex degree $d$ polynomials is
by definition the space of affine conjugacy classes of these polynomials.
An important subset of the parameter space is the so-called
\emph{connectedness locus} $\Mc_d$ consisting of classes of all
degree $d$ polynomials $P$, whose Julia sets $J(P)$ are connected.
For $d=2$, we obtain the famous \emph{Mandelbrot set} $\Mc_2$,
which can be identified with the set of complex numbers $c$ such that
$0$ does not escape to infinity under the iterations of the polynomial
$P_c(z)=z^2+c$.
The identification is based on the fact that every quadratic polynomial
is affinely conjugate to $P_c$ for some $c\in\C$ as well as a classical
theorem of Fatou and Julia.

\subsection{Combinatorial model of the Mandelbrot set}
A combinatorial model for $\Mc_2$ is due to Thurston \cite{thu85}.
It is constructed as follows.
Let $\uc$ be the unit circle in the plane of complex numbers, consisting
of all complex numbers of modulus one, and let $\si_2:\uc\to\uc$ be the angle-doubling
map $z\mapsto z^2$.
We will identify $\uc$ with $\R/\Z$ by means of the mapping taking an
\emph{angle} $\theta\in\R/\Z$ to the point $e^{2\pi i\theta}\in\uc$.
Under this identification, we have $\si_2(\theta)=2\theta$.
If the Julia set $J(P_c)$ is locally connected, then Thurston associates
a certain set $\lam_c$ of pairwise disjoint chords in the unit disk $\disk=\{z\in\C\,|\,|z|<1\}$
with the following property: the quotient space of the unit circle
$\uc/\lam_c$ obtained by identifying all pairs of points connected by chords in $\lam_c$
is homeomorphic to $J(P_c)$; moreover, the dynamics of $\si_2:\uc\to\uc$
descends to the quotient space, and the induced dynamics is topologically
conjugate to $P_c:J(P_c)\to J(P_c)$.

The set $\lam_c$ is called the \emph{geolamination}
(\emph{geodesic}, or \emph{geometric}, lamination) of $P_c$.
Thurston's geolaminations provide models for the topological dynamics of
quadratic polynomials with locally connected Julia sets.
It makes sense to consider limits of geolaminations $\lam_c$; these limits
(called \emph{limit quadratic geolaminations})
do not necessarily correspond to polynomials with locally connected Julia sets.
Chords belonging to a geolamination $\lam$ are called \emph{leaves} of $\lam$.
The main property that the leaves of a geolamination have is that
they are not \emph{linked}, i.e., they do not cross in $\disk$.

So far, this construction provides topological models for individual
quadratic polynomials --- not even for all of them, since there are polynomials
$P_c$ such that $J(P_c)$ is connected but not locally connected;
however, we need to model the space of \emph{all} polynomials $P_c$ with
connected Julia sets.
Metaphorically speaking, there are two parallel worlds: the ``analytic'' world of
complex polynomials and the ``combinatorial'' world of limit geolaminations.
Both worlds often come close to each other: whenever
we have a polynomial $P_c$ with locally connected $J(P_c)$, then we have the
corresponding geolamination $\lam_c$.
On the other hand, sometimes the two worlds diverge.
Still, a conjectural model for $\Mc_2$ can be built within the combinatorial world.

The idea is to take one particular leaf from every limit quadratic
geolamination $\lam$, namely, the leaf, called the \emph{minor} of $\lam$,
whose endpoints are the $\si_2$-images of the endpoints of a longest leaf of $\lam$.
The minors of all limit quadratic geolaminations form the so-called
\emph{quadratic minor lamination} $\mathrm{QML}$.
This is the geolamination that gives a conjectural model for the Mandelbrot set,
in the sense that the boundary of $\Mc_2$ is conjecturally homeomorphic to $\uc/\qml$.
The leaves of $\qml$ can be described without referring to limit
quadratic geolaminations.
To this end, let us first agree to denote by $|x-y|, x, y\in \uc$ the length of the shortest
circle arc with endpoints $x$ and $y$. Denote by $\ol{ab}$ the chord with endpoints $a$ and $b$.
Consider a chord $\ol{ab}$ with $|a-b|<1/3$.
Let $A$ be the shortest closed arc bounded by $a$ and $b$, and $S$ be
the convex hull of the set $\si_2^{-1}(A)$ in the plane.
The set $S$ is called the \emph{critical strip} of $\ell$.
A chord $\ell=\ol{ab}$ with endpoints $a$ and $b$ is a \emph{major}
if the following property holds:
for every positive integer $n$, the chord $\si_2^n(\ell)$ connecting
the points $\si_2^n(a)$ and $\si_2^n(b)$ is disjoint from the interior of $S$.
An alternative (and more straightforward) way of defining $\mathrm{QML}$ is saying
that $\qml$ is formed by $\si_2(\ell)$ for all majors $\ell$.
Note that the conjecture that the boundary of $\Mc_2$ is homeomorphic to $\uc/\qml$
is equivalent to the celebrated \emph{MLC conjecture} claiming that the
Mandelbrot set is locally connected.

We will write $\bd(X)$ for the boundary of a subset $X\subset\C$.
There is a continuous monotone mapping $\pi:\bd(\Mc_2)\to\uc/\qml$.
Recall that a continuous mapping from one continuum to another continuum
is \emph{monotone} if the \emph{fibers} (i.e., point preimages) are connected.
The set $\Mc_2$ is locally connected if and only if the fibers of $\pi$ are points,
hence, $\pi$ is the desired homeomorphism between $\bd(\Mc_2)$ and $\uc/\qml$
provided that the MLC conjecture holds.

The connectedness locus $\Mc_3$ in the parameter space of complex cubic
polynomials is a four-dimensional set
which is known to be non-locally connected \cite{la89}.
Thus, it is hopeless to look for a precise topological
model for the boundary of $\Mc_3$ as a quotient of a nice space like
the 3-sphere (any quotient space of a locally connected space is
locally connected!).
However, extensions of Thurston's results to the cubic case are possible
if, say, we study a rich enough subset of $\Mc_3$ instead of
the entire connectedness locus and if we allow for \emph{monotone} models
rather than precise ones.

In this paper, we study the space of cubic \emph{dendritic} polynomials.
These are polynomials with connected Julia sets, all of whose cycles are repelling.
Dendritic polynomials exhibit rich dynamics and have been actively studied before.
In particular, there is a
nice association, due to Kiwi \cite{kiwi97},
between dendritic polynomials and a certain class of geolaminations.

\subsection{Tagging dendritic cubic polynomials}
Similarly to the projection $\pi:\bd(\Mc_2)\to\uc/\qml$, we would like
to define a projection from the set of dendritic cubic polynomials
to a certain set of combinatorial objects.
The latter should be thought of as \emph{tags} of the dendritic polynomials.
The process of tagging is a two-step process.
Firstly, we associate every dendritic polynomial with the corresponding geolamination.
Secondly, we define a combinatorial tag of every ``dendritic'' geolamination.

The first step is essentially due to Jan Kiwi.
He showed in \cite{kiwi97} that, for every dendritic polynomial $P$ of degree $d$,
there is a monotone semi-conjugacy $\Psi_P$ between $P:J(P)\to J(P)$
and a certain quotient of $\si_d:\uc\to\uc$ represented by a geolamination $\lam_P$.
Here $\si_d$ is the $d$-tupling map $\theta\mapsto d\theta$ on the unit
circle, which descends to the quotient space $\uc/\lam_P$.
The corresponding induced continuous mapping $f:\uc/\lam_P\to\uc/\lam_P$ is called
the \emph{topological polynomial} associated with $P$.
As was already mentioned above, the quotient space $\uc/\lam_P$ is to be
understood as the quotient space of $\uc$ by a certain equivalence relation $\sim_P$.
By definition, $\sim_P$ is the minimal equivalence relation on $\uc$ with the
property that any two points connected by a leaf of $\lam_P$ are equivalent.
It turns out that, in the dendritic case, all classes of the equivalence relation
$\sim_P$ are finite.

Let us now discuss the second step, namely, the tagging of the geolaminations
$\lam_P$, or, equivalently, of the corresponding equivalence relations $\sim_P$.
We start again with the quadratic case.
Let $P_c(z)=z^2+c$ be a quadratic dendritic polynomial.
The corresponding parameter value $c$ is also called \emph{(quadratic) dendritic}.
Consider the $\sim_{P_c}$-equivalence class represented by the point
$\Psi_{P_c}(c)$ of $\uc/\sim_{P_c}$.
Let $G_c$ denote the convex hull of this class.
This is a convex polygon in the closed unit disk with finitely many vertices
on the unit circle.
This polygon may degenerate into a chord (if there are two vertices)
or even into a point (if there is just one vertex).
The fundamental results of Thurston imply, in particular, that $G_c$ and
$G_{c'}$ are either the same or disjoint, for all pairs $c$, $c'$ of
dendritic parameter values.
Moreover, the mapping $c\mapsto G_c$ is upper semicontinuous in a natural
sense (if a sequence of dendritic parameters $c_n$ converges to a
dendritic parameter $c$, then the limit set of the corresponding
convex sets $G_{c_n}$ is a subset of $G_c$). We call $G_c$ the
\emph{tag associated to $c$}.

Now, consider the union of all tags of quadratic dendritic
polynomials. This union is naturally partitioned into individual
tags (distinct tags are pairwise disjoint!). This defines its
quotient space. On the other hand, take the set of quadratic
dendritic parameters. Each such parameter $c$ maps to the polygon
$G_c$, i.e. to the tag associated to $c$. Thus, each quadratic
dendritic parameter maps to the corresponding point of the quotient
space of the union of all tags of quadratic dendritic polynomials
defined in the beginning of this paragraph. This provides for a
model of the set of quadratic dendritic polynomials (or their
parameters).

A major part of this paper is an extension of these results to the
cubic case. To explain our approach, we need a few definitions,
including some that can be useful in a more general setting.

Consider a dendritic polynomial $P$ of any degree.
We have the combinatorial objects $\lam_P$ and $\sim_P$ associated with $P$.
Given a point $z\in J(P)$, we
associate with it the convex hull $G_{P,z}$ of the $\sim_P$-equivalence
class represented by the point $\Psi_{P}(z)\in\uc/\sim_P$
(if $P$ is fixed, we may write $G_z$ instead of $G_{P,z}$).
The set $G_z$ is a convex polygon with finitely many vertices, a chord, or
a point; it should be viewed as a combinatorial object corresponding to $z$. For
any points $z\ne w\in J(P)$, the sets $G_z$ and $G_w$ either
coincide or are disjoint.

Let us now go back to cubic polynomials.
A \emph{critically marked} cubic polynomial is by definition a
triple $(P,\om_1,\om_2)$, where $P$ is a cubic polynomial with critical
points $\om_1$ and $\om_2$ such that $\om_1\ne\om_2$ unless $P$
has only one (double) critical point.
If $\om_1\ne \om_2$, then the triple $(P,\om_1,\om_2)$
and the triple $(P,\om_2,\om_1)$
are viewed as two distinct critically marked cubic polynomials.
Slightly abusing the notation, we will sometimes write $P$ for a critically
marked polynomial $(P,\om_1,\om_2)$, and then write $\om_i(P)$ instead of $\om_i$
to emphasize the dependence on $P$.
Let $\mathcal{MD}_3$ be the space of all critically marked cubic dendritic
polynomials. The \emph{co-critical} point
associated to a critical point $\om_i=\om_i(P)$ of a cubic polynomial $P$ is
the only point $\om_i^*$ with $P(\om_i^*)=P(\om_i)$ and $\om_i^*\ne\om_i$
unless $P$ has a double critical point $\om$ in which case $\om^*=\om$.
Then, with every marked dendritic polynomial $P$, we associate the corresponding
\emph{mixed tag}
$$
\ta(P)=G_{\om_1(P)^*}\times G_{P(\om_2(P))}\subset
\overline{\disk}\times\overline{\disk}.
$$

Let $\ta(\md_3)^+$ be the union of the sets $\ta(P)$ over all $P\in\md_3$.
It turns out that the mixed tags $\ta(P)$ form a partition of $\ta(\md_3)^+$
and generate the corresponding quotient space of $\ta(\md_3)^+$
denoted by $\cml$ (for \emph{cubic mixed lamination}). Moreover, we prove that
$\ta:\md_3\to \cml$ is continuous and thus $\cml$ can serve as a
combinatorial model for $\md_3$. All this is summarized below in our Main Theorem.

\begin{thmM}
Mixed tags of critically marked polynomials from $\mathcal{MD}_3$ are disjoint
or coincide. The map $\ta: \mathcal{MD}_3\to \cml$ is continuous.
\end{thmM}

Thus, there is a continuous mapping from the space of marked cubic
dendritic polynomials to the model space of their tags defined
through (geo)\-la\-mi\-na\-tions associated with marked polynomials
from $\mathcal{MD}_3$. This can be viewed as a partial
generalization of \cite{thu85} to cubic polynomials.

\subsection{Previous work}
Branner and Hubbard \cite{BrHu} initiated the study of $\Mc_3$, and studied 
the complement of this set in the full parameter space of cubic polynomials.
The complement is foliated by so-called \emph{stretching rays}
that are in a sense analogous to external rays of the Mandelbrot set.
The combinatorics of $\Mc_3$ is closely related to landing patterns
of stretching rays. However, we do not explore this connection here.
A significant complication is caused 
by the fact that there are
many non-landing stretching rays.
Landing properties of stretching rays in the parameter space of real
polynomials have been studied by Komori and Nakane \cite{kona}.

Lavaurs \cite{la89} proved that
$\Mc_3$ is not locally connected. Epstein and Yampolsky \cite{EY99}
proved that the bifurcation locus in the space of real cubic
polynomials is not locally connected either. This makes the problem
of defining a 
combinatorial model of $\Mc_3$ very
delicate. Buff and Henriksen \cite{BH01} presented copies of
quadratic Julia sets, including Julia sets that are not locally
connected, in slices of $\Mc_3$. In his thesis, D. Faught
\cite{Fau92} considered the slice $\mathcal A$ of $\Mc_3$ consisting of
polynomials with a fixed critical point and showed that $\mathcal A$ contains
countably many homeomorphic copies of $\Mc_2$ and is locally
connected everywhere else. 
P. Roesch \cite{R06} filled the gaps in Faught's arguments and
generalized Faught's results to higher degrees.
Milnor \cite{M} gave a classification of hyperbolic components in $\Mc_d$;
however, this description does not involve combinatorial tags.
Schleicher \cite{sch04} constructed a geolamination modeling the
space of \emph{unicritical} cubic polynomials, i.e., cubic polynomials with
a multiple critical point.
We have heard of an unpublished old work of D. Ahmadi and M. Rees,
in which cubic geolaminations were studied, however we have not seen
it. 

\subsection{Overview of the method}

Thuston's tools used in the construction of $\qml$
do not generalize to the cubic case.
These tools are based on the Central Strip Lemma stated in Section \ref{ss:motiv},
and include the No Wandering Triangles Theorem
(also stated in Section \ref{ss:motiv}).
A straightforward extension of the Central Strip Lemma as well as that of the No
Wandering Triangles Theorem to the cubic case fail, e.g., cubic geolaminations
may have wandering triangles, cf. \cite{bl02}.
Thus, one needs a different set of combinatorial tools.
Such tools are developed in this paper --- they are called \emph{smart criticality}.
Smart criticality works for geolaminations of any degree.

Given a geolamination $\lam$, define \emph{gaps} of $\lam$ as closure of components
of $\cdisk\sm \lam^+$ where $\lam^+\subset \cdisk$ is the union of all
leaves of $\lam$.
The statement about the quadratic laminations we are trying to generalize
is the following: if the minors of two quadratic geolaminations intersect
in $\disk$, then they coincide.
Although minors can also be defined for higher degree laminations, they are
not the right objects to consider.
For a quadratic geolamination $\lam$, instead of its non-degenerate minor $\ol m$, we can
consider the quadrilateral, whose vertices are the four $\si_2$-preimages
of the endpoints of $\ol m$.
Such a quadrilateral is called a \emph{critical quadrilateral}.
The critical quadrilateral of a quadratic geolamination $\lam$ lies in
some gap of $\lam$ or, if $\ol m$ is a point, coincides with a leaf of $\lam$.
Similarly, for a degree $d$ invariant geolamination $\lam$, we can define
critical quadrilaterals as quadrilaterals (possibly degenerate) lying in gaps
or leaves of $\lam$, whose opposite vertices have the same $\si_d$-images.
These critical quadrilaterals will play the role of minors and will be
used to tag higher degree geolaminations.

The method of smart criticality helps to verify that, under suitable assumptions,
two linked leaves $\ell_1$, $\ell_2$ of \emph{different} geolaminations
have linked images $\si_d^n(\ell_1)$, $\si_d^n(\ell_2)$, for \emph{all} $n$.
One possible reason, for which $\si_d(\ell_1)$, $\si_d(\ell_2)$ may be linked,
is the following: $\ell_1$ and $\ell_2$ are disjoint from a full collection
of critical chords (here
a $\si_d$-\emph{critical chord} is a chord of $\disk$, whose endpoints map to
the same point under $\si_d$, and
a \emph{full collection of critical chords} is a collection of $d-1$ critical chords
without loops).
To prove that $\si_d^n(\ell_1)$, $\si_d^n(\ell_2)$ are 
linked for all $n$,
we will choose, for every $n$, a different full collection of critical chords ---
this is the meaning of ``smart''.

Smart criticality can be implemented in the following situation.
Let $\lam_1$ and $\lam_2$ be two geolaminations.
Suppose that we can choose critical quadrilaterals in $\lam_1$ and $\lam_2$
so that the corresponding quadrilaterals of different geolaminations either
have alternating vertices, or share a diagonal.
In this case, we say that $\lam_1$ and $\lam_2$ are \emph{linked or essentially equal}.
In fact, being linked or essentially equal is slightly more general than
the property just stated; the precise definition is Definition \ref{d:qclink1}.
Suppose now that $\lam_1$ and $\lam_2$ correspond to dendritic polynomials.
Smart criticality implies that, if $\lam_1$ and $\lam_2$ are
linked or essentially equal, then they must coincide.
Together with some purely combinatorial (and non-dynamical) considerations,
this translates into the following statement:
if the tags of $\lam_1$ and $\lam_2$ are non-disjoint, then $\lam_1=\lam_2$.
Basically, this is all we need in order to prove the Main Theorem.

Our main tools (smart criticality) are developed for geolaminations of any degree.
However, the Main Theorem is confined with cubic polynomials and
cubic geolaminations.
The reason is that the purely combinatorial and non-dynamical
considerations that help to translate non-disjointness of tags
into the linkage of geolaminations are much 
more involved in the higher degree case. Thus, even though
we believe that the Main Theorem generalizes to all degrees,
a lot of details will have to be worked out and a
careful proof would require a significant additional space and time.

\subsection{Organization of the paper}
In Section \ref{s:basicdef}, we discuss general properties of geolaminations
as well as specific classes of geolaminations, e.g., dendritic geolaminations.
We also introduce combinatorial objects (qc-portraits) that serve
as combinatorial tags of geolaminations.
In Section \ref{s:acclam}, we study so-called \emph{accordions}.
These are geometric objects formed by crossing leaves of different geolaminations.
Smart criticality yields that accordions of linked or essentially equal
geolaminations behave much like gaps of a single geolamination.
This is established in Section \ref{s:qcrit}, where the method of smart
criticality is developed.
Finally, in Section \ref{l:appli}, we will prove the Main Theorem.

\section{Geolaminations and their properties}\label{s:basicdef}

In this section, we give basic definitions, list some known results
on geolaminations, and establish some new facts about them.

\subsection{Basic definitions}\label{ss:bd}
For a collection $\rc$ of chords of $\disk$ set $\bigcup\rc=\rc^+$.
A \emph{geolamination} is a collection $\lam$ of (perhaps
degenerate) chords of $\cdisk$ called \emph{leaves} which are
pairwise disjoint in $\disk$ such that
$\lam^+=\bigcup_{\ell\in\lam}\ell$ is closed, and all points of
$\uc$ are elements of $\lam$.  We linearly extend $\si_d$ over
leaves of $\lam$; clearly, this extension is continuous and
well-defined.
We define \emph{gaps} of $\lam$ as the closures of the components
of $\disk\sm\lam^+$.

\subsubsection{Sibling invariant geolaminations}\label{sss:sibgeo}

Let us introduce the notion of a (sibling) \emph{$\si_d$-invariant}
geolamination which is a slight modification of an invariant
geolamination introduced by Thurston \cite{thu85}.

\begin{dfn}[Invariant geolaminations \cite{bmov13}]\label{d:sibli}
A geolamination $\lam$ is (sibling) \emph{$\si_d$-invariant}
provided that:
\begin{enumerate}
\item for each $\ell\in\lam$, we have $\si_d(\ell)\in\lam$,
\item \label{2}for each $\ell\in\lam$ there exists $\ell^*\in\lam$ so that $\si_d(\ell^*)=\ell$.
\item \label{3} for each $\ell\in\lam$ such that $\si_d(\ell_1)$
    is a non-degenerate leaf, there exist $d$ \textbf{pairwise disjoint}
 leaves $\ell_1$, $\dots$, $\ell_d$ in $\lam$ such that
 $\ell_1=\ell$ and
 $\si_d(\ell_i)=\si_d(\ell)$ for all
  $i=2$, $\dots$, $d$.
\end{enumerate}
\end{dfn}

We call the leaf $\ell^*$ in (\ref{2}) a \emph{pullback} of $\ell$
and the leaves $\ell_2$, $\dots$, $\ell_d$ in (\ref{3})
\emph{siblings} of $\ell=\ell_1$. In a broad sense a \emph{sibling
of $\ell$} is a leaf with the same image but distinct from $\ell$.
Definition~\ref{d:sibli} is slightly more restrictive than
Thurston's definition of an invariant geolamination. By
\cite{bmov13}, a $\si_d$-invariant geolamination $\lam$ is
\emph{invariant} in the sense of Thurston \cite{thu85} and, in
particular, \emph{gap invariant}: if $G$ is a gap of $\lam$ and $H$
is the convex hull of $\si_d(G\cap \uc)$, then $H$ is a point, a
leaf of $\lam$, or a gap of $\lam$, and in the latter case, the map
$\si_d|_{\bd(G)}:\bd(G)\to \bd(H)$ of the boundary of $G$ onto the
boundary of $H$ is a positively oriented composition of a monotone
map and a covering map. From now on by \emph{$(\si_d$-$)$invariant}
geolaminations we mean sibling $\si_d$-invariant geolaminations and
consider \emph{only} such invariant geolaminations.

\begin{thm}[Theorem 3.21 \cite{bmov13}]\label{t:sibliclos}
The family of sets $\lam^+$ of all invariant geolaminations $\lam$
is closed in the Hausdorff metric. In particular, this family is compact.
\end{thm}

Clearly, $\lam_i^+\to \lam^+$ (understood as convergence of compact
subsets of $\cdisk$) implies that the collections of chords $\lam_i$
converge to the collection of chords $\lam$ (i.e., each leaf of
$\lam$ is the limit of a sequence of leaves from $\lam_i$, and each
converging sequence of leaves of $\lam_i$ converges to a leaf of
$\lam$). Thus, from now on we will write $\lam_i\to \lam$ if
$\lam^+_i\to \lam^+$ in the Hausdorff metric.

Two \emph{distinct} chords of $\disk$ are \emph{linked} if
they intersect in $\disk$ (we will also say that these chords
\emph{cross each other}). A gap $G$ is called \emph{infinite
$($finite, uncountable$)$} if $G\cap \uc$ is infinite (finite,
uncountable). Uncountable gaps are also called \emph{Fatou} gaps.
For a closed convex set $H\subset \C$, straight segments from $\bd(H)$ are called
\emph{edges} of $H$. The degree of a gap or leaf $G$ is defined as follows.
If $\si_d(G)$ is degenerate then the degree of $G$ is the cardinality
of $G\cap \uc$. Suppose now that $\si_d(G)$ is not a point.
Consider $\si_d|_{\bd(G)}$. Then the degree of $G$ equals the number of
components in the preimage of a point $z\in \si_d(\bd(G))$ under the map
$\si_d|_{\bd(G)}$.

\begin{dfn}\label{d:cristuff} We say
that $\ell$ is a \emph{chord of a geolamination $\lam$} if $\ell$ is
a chord of $\ol\disk$ unlinked with all leaves of $\lam$. A
\emph{critical chord $($leaf$)$} $\ol{ab}$ of $\lam$ is a chord
(leaf) of $\lam$ such that $\si_d(a)=\si_d(b)$. A gap is
\emph{all-critical} if all its edges are critical. An all-critical
gap or a critical leaf is called an \emph{all-critical set}. A gap
$G$ is said to be \emph{critical} if the degree of $G$ is greater
than one. A \emph{critical set} is either a critical leaf or a
critical gap.
\end{dfn}

By Thurston \cite{thu85}, there is a canonical \emph{barycentric}
extension of the map $\si_d$ to the entire closed disk $\cdisk$.
First $\si_d$ is extended linearly over all leaves of an invariant
geolamination $\lam$, and then piecewise linearly over the interiors
of all gaps of $\lam$, using the barycentric subdivision.
When talking about $\si_d$ on $\cdisk$, we always have some invariant
geolamination in mind and mean Thurston's barycentric extension
described above.

\subsubsection{Laminations as equivalence relations}\label{sss:laeqre} A lot of
geolaminations naturally appear in the context of invariant
equivalence relations on $\uc$ (\emph{laminations}) satisfying
special conditions.

\begin{dfn}[Laminations]\label{d:lam}
An equivalence relation $\sim$ on the unit circle $\uc$ is called a
\emph{lamination} if either $\uc$ is one $\sim$-class (such
laminations are called \emph{degenerate}),
or the following holds:

\noindent (E1) the graph of $\sim$ is a closed subset of $\uc \times
\uc$;

\noindent (E2) the convex hulls of distinct equivalence classes are
disjoint;

\noindent (E3) each equivalence class of $\sim$ is finite.
\end{dfn}

\begin{dfn}[Laminations and dynamics]\label{d:si-inv-lam}
An equivalence relation $\sim$ is called ($\si_d$-){\em invariant}
if:

\noindent (D1) $\sim$ is {\em forward invariant}: for a $\sim$-class
$\g$, the set $\si_d(\g)$ is a $\sim$-class;

\noindent (D2) for any $\sim$-class $\g$, the map $\si_d: \g\to
\si_d(\g)$ extends to $\uc$ as an orientation preserving covering
map such that $\g$ is the full preimage of $\si_d(\g)$ under this
covering map.
\end{dfn}

For an invariant lamination $\sim$ consider the \emph{topological
Julia set} $\uc/\hspace{-5pt}\sim\,=J_\sim$ and the
\emph{topological polynomial} $f_\sim:J_\sim\to J_\sim$ induced by
$\si_d$. The quotient map $\pi_\sim:\uc\to
\uc/\hspace{-5pt}\sim=J_\sim$ semi-conjugates $\si_d$ with
$f_\sim|_{J_\sim}$. A lamination $\sim$ admits a \emph{canonical
extension over $\C$}: nontrivial classes of this extension are convex
hulls of classes of $\sim$. By Moore's Theorem, the quotient space
$\C/\hspace{-5pt}\sim$ is homeomorphic to $\C$. The quotient map
$\pi_\sim:\uc\to \uc/\hspace{-5pt}\sim$ extends to the plane with
the only non-trivial point-preimages (\emph{fibers}) being the
convex hulls of non-degenerate $\sim$-classes.
With any fixed identification between $\C/\sim$ and $\C$, one can extend
$f_\sim$ to a branched-covering map $f_\sim:\C\to \C$ of degree $d$
called a \emph{topological polynomial} too. The complement $K_\sim$
of the unique unbounded component $U_\infty(J_\sim)$ of $\C\sm
J_\sim$ is called the \emph{filled topological Julia set}. Define the \emph{canonical geolamination
$\lam_\sim$ generated by $\sim$} as the collection of edges of
convex hulls of all $\sim$-classes and all points of $\uc$. By
\cite{bmov13}, the geolamination $\lam_\sim$ is $\si_d$-invariant.

\subsubsection{Other useful notions}\label{sss:oun}
Considering objects related to (geo)laminations, we do not have
to fix these (geo)laminations.

\begin{dfn}\label{d:restuff}
By a \emph{periodic gap or leaf}, we mean a gap or a leaf $G$, for which
there exists the least number $n$ (called the \emph{period} of $G$)
such that $\si_d^n(G)=G$. Then we call the map $\si_d^n:G\to G$ the
\emph{remap}. An edge (vertex) of $G$ on which the remap is identity is said to be
\emph{refixed}.
\end{dfn}

Given two points $a$, $b\in \uc$ we denote by $(a,b)$ the
\emph{positively oriented} arc from $a$ to $b$ (i.e., moving from
$a$ to be $b$ within $(a, b)$ takes place in the counterclockwise
direction). For a closed set $G'\subset \uc$,
we call components of $\uc\sm G'$ \emph{holes}. If $\ell=\ol{ab}$ is
an edge of $G=\ch(G')$, then we let
$H_G(\ell)$ denote the component
of $\uc\sm\{a, b\}$ disjoint from $G'$ and call it the hole of $G$
\emph{behind $\ell$} (it is only unique if $G'$ contains at least three points). The \emph{relative interior} of a gap is its
interior in the plane; the \emph{relative interior} of a segment is
the segment minus its endpoints.

\begin{dfn}\label{d:lamset}
If $A\subset \uc$ is a closed set such that all the sets
$\ch(\si_d^i(A))$ are pairwise disjoint, then $A$ is called
\emph{wandering}. If there exists $n\ge 1$ such that all the sets
$\ch(\si_d^i(A)), i=0, \dots, n-1$ have pairwise disjoint relative
interiors while $\si_d^n(A)=A$, then $A$ is called \emph{periodic}
of period $n$. If there exists $m>0$ such that all $\ch(\si_d^i(A)),
0\le i\le m+n-1$ have pairwise disjoint relative interiors and
$\si_d^m(A)$ is periodic of period $n$, then we call $A$
\emph{preperiodic} of period $n$ and preperiod $m$. If $A$ is
wandering, periodic or preperiodic, and for every $i\ge 0$ and every
hole $(a, b)$ of $\si_d^i(A)$ either $\si_d(a)=\si_d(b)$, or the
positively oriented arc $(\si_d(a), \si_d(b))$ is a hole of
$\si_d^{i+1}(A)$, then we call $A$ (and $\ch(A)$) a
\emph{{\rm(}$\si_d${\rm )}-laminational set}; we call $\ch(A)$
\emph{finite} if $A$ is finite. A {\em {\rm(}$\si_d$-{\rm)}stand
alone gap} is defined as a laminational set with non-empty interior.
\end{dfn}

Denote by $<$ the \emph{positive} (counterclockwise) circular order
on $\uc=\mathbb R/\mathbb Z$ induced by the usual order of $\mathbb
R$. Note that this order is only meaningful for sets of cardinality
at least three. For example, we say that $x<y<z$ provided that
moving from $x$ in the positive direction along $\uc$ we meet $y$
before meeting $z$.

\begin{dfn}[Order preserving] Let $X\subset \uc$ be a set with at
least three points. We call $\si_d$ \emph{order preserving on $X$}
if $\si_d|_X$ is one-to-one and, for every triple $x$, $y$, $z\in X$
with $x<y<z$, we have $\si_d(x)<\si_d(y)<\si_d(z)$.
\end{dfn}

\subsection{General properties of invariant geolaminations}\label{ss:gpgeo}

\begin{lem}[Lemma 3.7 \cite{bmov13}] \label{l:3.7}
If $\ol{ab}$ and $\ol{ac}$ are two leaves of an invariant geolamination $\lam$
such that $\si_d(a), \si_d(b)$ and $\si_d(c)$ are all distinct
points, then the order among points $a$, $b$, $c$ is preserved under $\si_d$.
\end{lem}

We prove a few corollaries of Lemma~\ref{l:3.7}

\begin{lem}\label{l:sameperiod1}
If $\lam$ is an invariant geolamination, $\ell=\ol{ab}$ is a leaf of $\lam$, and
the point $a$ is periodic, then $b$ is $($pre$)$periodic of the same period.
\end{lem}

\begin{proof}
Assume that $a$ is of period $n$ but $b$ is not $\si_d^n$-fixed. Then, by Lemma
\ref{l:3.7}, either the circular order among the points $b_i=\si_d^{ni}(b)$ is
the same as the order of subscripts or $b_i=b_{i+1}$ for some
$i$. In the former case $b_i$ converge to some limit point, a
contradiction with the expansion property of $\si^n_d$. Hence for
some (minimal) $i$ we have $b_i=b_{i+1}$. It follows that the period $m$
of $b_i$ cannot be less than $n$ as otherwise we can consider
$\si_d^m$ which fixes $b_i$ and does not fix $a$ yielding the same contradiction
with Lemma \ref{l:3.7}.
\end{proof}

We will need the following elementary lemma.

\begin{lem}\label{l:concat}
If $x\in \uc$ and the chords
$\ol{\si^i_d(x)\si_d^{i+1}(x)}$, $i=0$, $1$, $\dots$ are pairwise
unlinked then $x$ $($and hence the leaf $\ol{x\si_d(x)}=\ell)$ is
$($pre$)$periodic.
\end{lem}

\begin{proof}
The sequence of leaves from the lemma is the $\si_d$-orbit of $\ell$,
in which consecutive images are concatenated and no two leaves are
linked. If, for some $i$, the leaf
$\ol{\si^i_d(x)\si_d^{i+1}(x)}=\si^i_d(\ell)$ is critical, then
$\si^{i+1}_d(\ell)=\{\si_d^{i+1}(x)\}$ is a $\si_d$-fixed point,
which proves the claim in this case. Assume now that $\ell$ is not
(pre)critical. If $x$ is not (pre)periodic, then, by topological
considerations, leaves $\si_d^n(\ell)$ must converge to a limit leaf
or point. Clearly, this limit set is $\si_d$-invariant. However, $\si_d$
is expanding, a contradiction.
\end{proof}

Lemma~\ref{l:concat} easily implies Lemma~\ref{l:conconv}.

\begin{lem}\label{l:conconv}
Let $\lam$ be a geolamination. Then the following holds.

\begin{enumerate}

\item
If $\ell$ is a leaf of $\lam$ and, for some $n>0$, the leaf
$\si_d^n(\ell)$ is concatenated to $\ell$, then $\ell$ is
$($pre$)$periodic.

\item If $\ell$ has a (pre)periodic endpoint, then $\ell$ is
$($pre$)$periodic.

\item If two leaves $\ell_1$, $\ell_2$ from geolaminations $\lam_1$,
$\lam_2$ share the same $($pre$)$\-pe\-rio\-dic endpoint, then they
are $($pre$)$periodic with the same eventual period of their
endpoints.

\end{enumerate}

\end{lem}

\begin{proof}
Let $\ell=\ol{uv}$.
First, assume that $\si_d^n(u)=u$. Then (1) follows from Lemma
\ref{l:sameperiod1}. Second, assume that $\si^n_d(u)=v$. Then (1)
follows from Lemma~\ref{l:concat}. Statements (2) and (3) follow
from (1) and Lemma~\ref{l:sameperiod1}. \end{proof}

A similar conclusion can be made for edges of periodic gaps.

\begin{lem}\label{l:cripe}
Any edge of a periodic gap is $($pre$)$pe\-ri\-odic or
$($pre$)$\-cri\-ti\-cal.
\end{lem}

\begin{proof}
Let $G$ be a fixed gap and $\ell$ be a non-(pre)critical edge of it.
The length $s_n$ of the hole $H_G(\si_d^n(\ell))$ of
$G$ behind the leaf $\si_d^n(\ell)$ grows with $n$ as long as $s_n$
stays sufficiently small (it is easy to see that the correct bound
on $s_n$ is that $s_n<\frac{1}{d+1}$). Hence the
sequence $\{s_i\}$ will contain infinitely many numbers greater than
or equal to $\frac 1{d+1}$. A contradiction with the fact that there are
only finitely many distinct holes of $G$ of length $\frac{1}{d+1}$
or bigger.
\end{proof}

Given $v\in \uc$, let $E(v)$ be the closure of the set $\{u\,|\, \ol{uv}\in\lam\}$.

\begin{lem}\label{l:cones}
If $v$ is not $($pre$)$periodic, then $E(v)$ is at most finite. If
$v$ is $($pre$)$periodic, then $E(v)$ is at most countable.
\end{lem}

\begin{proof}
The first claim is proven in \cite[Lemma 4.7]{bmov13}. The second
claim follows from Lemma~\ref{l:conconv}.
\end{proof}

Properties of individual wandering polygons were studied in
\cite{kiw02}; properties of collections of wandering polygons were
studied in \cite{bl02}; their existence was established in
\cite{bo08}. The most detailed results on wandering polygons and
their collections are due to Childers \cite{chi07}.

Let us describe the entire $\si_d$-orbit of a finite periodic
laminational set.

\begin{prop}\label{p:forconcat}
Let $T$ be a $\si_d$-periodic finite laminational set and $X$ be the
union of the forward images of $T$. Then, for every connected
component $R$ of $X$, there is an $m$-tuple of points
$a_0<a_1<\dots<a_{m-1}<a_m=a_0$ in $\uc$ such that $R$ consists of
eventual images of $T$ containing $\ol{a_ia_{i+1}}$ for $i=0, \dots,
m-1$. If $m>1$, then the remap of $R$ is a combinatorial rotation
sending $a_i$ to $a_{i+1}$.
\end{prop}

Note that the case $m=1$ is possible. In this case, $R$ consists of
several images of $T$ sharing a common vertex $a_0$, there
is a natural cyclic order among the images of $T$, and the remap of
$R$ is a cyclic permutation of these images, not necessarily a
combinatorial rotation.

\begin{proof}
Set $T_k=\si_d^k(T)$. Let $k$ be the smallest positive integer such
that $T_k$ intersects $T_0$; we may suppose that $T_k\ne T_0$. There
is a vertex $a_0$ of $T_0$ such that $a_1=\si_d^k(a_0)$ is also a
vertex of $T_0$. Clearly, both $a_1$ and $a_2=\si_d^k(a_1)$ are
vertices of $T_k$. Set $a_i=\si_d^{ki}(a_0)$. Then we have $a_m=a_0$
for some minimal $m>0$. Let $Q$ be the convex hull of the points $a_0$,
$\dots$, $a_{m-1}$. This is a convex polygon, or a chord, or a
point. If $m>1$, then $a_i$ and $a_{i+1}$ are the endpoints of the
same edge of $Q$ (otherwise some edges of the polygons $T_{ki}$ would cross in
$\disk$). Set $R=\cup_{i=0}^{m-1}T_{ki}$. If $m=1$, then the sets
$T_{ki}$ share the vertex $a_0$. If $m>1$, then
every chord $\ol{a_ia_{i+1}}$ is an edge of $T_{ki}$ shared with $Q$, sets
$T_{ki}$ are disjoint from the interior of $Q$, and the remap
$\si_d^k$ of $R$ is a combinatorial rotation acting transitively on
the vertices of $Q$.

To prove that $R$ is disjoint from $R_j=\si_d^j(R)$ for $j<k$
suppose that $R_j$ intersects some $T_{ki}$. Note that the map
$\si_d^k$ fixes both $R$ and $R_j$. It follows that $R_j$ intersects
all $T_{ki}$, hence contains $Q$, a contradiction.
\end{proof}

It is well-known \cite{kiw02} that any infinite gap $G$ of a
geolamination $\lam$ is (pre)periodic. By a \emph{vertex} of a gap or leaf $G$ we
mean any point of $G\cap\uc$.

\begin{lem}\label{l:perinf}
Let $G$ be a periodic gap of period $n$ and set $K=\bd(G)$. Then
$\si_d^n|_K$ is 
the composition of a covering map and a monotone map of
$K$. If $\si_d^n|_K$ is of degree one, then either {\rm (1)} or {\rm
(2)} holds.

\begin{enumerate}
\item The gap $G$ has countably many vertices, only finitely many
of which are periodic. All non-periodic edges of $G$ are $($pre$)$critical.

\item The map $\si_d^n|_K$ is monotonically semiconjugate to an irrational
circle rotation so that each fiber of this semiconjugacy is a finite
concatenation of $($pre$)$critical edges of $G$.
\end{enumerate}

\end{lem}

\begin{proof}
We will prove only the very last claim. Denote by $\ph$ the
semiconjugacy from (2). Let $T\subset K$ be a fiber of $\ph$. By
Lemma~\ref{l:cripe} all edges of $G$ are (pre)critical. Hence if $T$
contains infinitely many edges, then the forward images of $T$ will
hit critical leaves of $\si_d^n$ infinitely many times as $T$ cannot
collapse under a finite power of $\si_d^n$. This would imply that an
irrational circle rotation has periodic points, a contradiction.
\end{proof}

Lemma~\ref{l:perinf} implies Corollary~\ref{c:noinf}.

\begin{cor}\label{c:noinf}
Suppose that $G$ is a periodic gap of a geolamination $\lam$, whose
remap has degree one. Then at most countably many pairwise unlinked
leaves of other geolaminations can be located inside $G$.
\end{cor}

We say that a chord is located \emph{inside} $G$ if it is a subset
of $G$ and intersects the interior of $G$.

\begin{proof}
Any chord located inside $G$ has its endpoints at vertices of $G$.
Since in case (1) of Lemma~\ref{l:perinf} there are countably many
vertices of $G$, we may assume that case (2) of Lemma~\ref{l:perinf}
holds. Applying the semiconjugacy $\varphi$ from this lemma we see
that if a leaf $\ell$ is located in $G$ and its endpoints do not map
to the same point by $\varphi$, then $\ell$ will eventually cross
itself. If there are uncountably many leaves of geolaminations inside $G$,
then among them there must exist a leaf $\ell$ with endpoints in
distinct fibers of $\ph$. By the above some forward images of $\ell$
cross each other, a contradiction.
\end{proof}

\subsection{Geolaminations with qc-portraits}\label{ss:qc-por}

Here we define geo\-la\-mi\-na\-tions with quadratically critical
(qc-)portraits and discuss linked or essentially coinciding
geolaminations with qc-portraits. First we motivate our approach.

Thurston defines the \emph{minor} $m$ of a $\si_2$-invariant
lamination $\lam$ as the image of a longest leaf $M$ of $\lam$.
Any longest leaf of $\lam$ is said to be a \emph{major} of $\lam$. If $m$ is non-degenerate,
$\lam$ has two disjoint majors which both map to $m$; if $m$ is degenerate,
$\lam$ has a unique major which is a critical leaf. In the quadratic case
the majors are uniquely determined by the minor. Even though in the
cubic case one could define majors and minors similarly, unlike in
the quadratic case these ``minors'' do not uniquely determine the
corresponding majors. The simplest way to see that is to consider
distinct pairs of critical leaves with the same images. One can
choose two all-critical triangles with so-called \emph{aperiodic
kneadings} as defined by Kiwi in \cite{kiwi97}. By \cite{kiwi97},
this would imply that any choice of two disjoint critical leaves,
one from either triangle, will give rise to the corresponding
geolamination; clearly, these two geolaminations are very different
even though they have the same images of their critical leaves,
i.e., the same minors. Thus, in the cubic case we should be concerned with critical
sets, not only their images.

We study how ordered collections of critical sets of geolaminations
are located with respect to each other. The fact that critical sets
may have different degrees complicates such study. So, it is natural
to adjust our geolaminations to make sure that the associated critical
sets of two geolaminations are of the same type.

\begin{dfn}\label{d:qll}
A (generalized) \emph{critical quadrilateral} $Q$  is the circularly
ordered 4-tuple $[a_0,a_1,a_2,a_3]$ of marked points $a_0\le a_1\le
a_2\le a_3\le a_0$ in $\uc$ so that $\ol{a_0a_2}$ and $\ol{a_1a_3}$ are
critical chords (called \emph{spikes}); here critical
quadrilaterals $[a_0,a_1,a_2,a_3]$, $[a_1,a_2,a_3,a_0]$, $[a_2,a_3,a_0,a_1]$
and $[a_3,a_0,a_1,a_2]$ are viewed as equal.
\end{dfn}

We want to comment upon our notation. By $(X_1, \dots, X_k)$, we
always mean a $k$-tuple, i.e., an \emph{ordered} collection of
elements $X_1, \dots, X_k$. On the other hand, by $\{X_1, \dots,
X_k\}$ we mean a collection of elements $X_1, \dots, X_k$ with no
fixed order. Since, for critical \ql s, we need to emphasize the
\emph{circular} order among its vertices, we choose the notation $[a_0,a_1,
a_2, a_3]$ distinct from either of the two just described notations.

For brevity, we will often use the expression ``critical \ql'' when
talking about the convex hull of a critical \ql{}. Clearly, if all
vertices of a critical \ql{} are distinct or if its convex hull is a critical leaf, then
the \ql{} is uniquely defined by its convex hull. However, if the convex
hull of a critical \ql{} is a triangle, this is no longer true. Indeed,
let $T=\ch(a, b, c)$ be an all-critical triangle. Then $[a,a,b,c]$ is a
critical \ql, but so are $[a,b,b,c]$ and $[a,b,c,c]$.

A \emph{collapsing quadrilateral} is a critical quadrilateral, whose
$\si_d$-image is a leaf. A critical quadrilateral $Q$ has two intersecting
spikes and is a collapsing quadrilateral, a critical leaf, an
all-critical triangle, or an all-critical quadrilateral. If all its vertices
are pairwise distinct, we call $Q$ \emph{non-degenerate}, otherwise
$Q$ is called \emph{degenerate}. Vertices $a_0$ and $a_2$ ($a_1$ and
$a_3$) are called \emph{opposite}. Considering geolaminations, all of
whose critical sets are critical \ql s, is not very restrictive: we
can ``tune'' a given geolamination by inserting new leaves into its
critical sets in order to construct a new geolamination with all
critical sets being critical \ql s.

\begin{lem}\label{l:qls}
The family of all critical \ql s is closed. The family of all
critical \ql s that are critical sets of geolaminations is closed
too.
\end{lem}

\begin{proof}
The first claim is trivial. The second one follows from
Theorem~\ref{t:sibliclos} and the fact that if $\lam_i\to \lam$,
then the critical \ql s of geolaminations $\lam_i$ converge to
critical \ql s that are critical sets of $\lam$.
\end{proof}

In the quadratic case we have less variety of critical \ql s: only
collapsing \ql s and critical leaves. As mentioned above, each
quadratic invariant geolamination $\lam$ either already has a
critical \ql, or can be tuned to have one. The latter can be done in
several ways if $\lam$ has a finite critical set (on which $\si_2$
acts two-to-one). If however $\lam$ does not have a finite critical
set, then its critical set must be a periodic Fatou gap $U$ of
degree two. It follows from \cite{thu85} that it has a unique
refixed edge $M$; then one can tune $\lam$ by inserting into $U$ the \ql{}
which is the convex hull of $M$ and its sibling.

Thurston's parameterization \cite{thu85} can be viewed as
associating to every geolamination $\lam$ with critical \ql{} $Q$
its minor $m$. It is easy to see that $m$ is the $\si_2$-image of $Q$
and that $Q$
is the full $\si_2$-preimage of $m$. We would like
to translate some crucial results of Thurston's into the language of
critical \ql s of quadratic geolaminations. To this end, observe, that,
by the above, two minors cross if and only if their full pullbacks
(which are collapsing \ql s coinciding with convex hulls of pairs of
majors) have a rather specific mutual location: their vertices
alternate on the circle. A major result of Thurston's from
\cite{thu85} is that \emph{minors of different quadratic
geolaminations are unlinked}; in the language of critical \ql s this
can be restated as follows: \emph{critical \ql s of distinct
quadratic geolaminations cannot have vertices which alternate on the
circle}. All this motivates Definition~\ref{d:strolin}.

\begin{dfn}\label{d:strolin}
Let $A$ and $B$ be two \ql s. Say that $A$ and $B$ are
\emph{strongly linked} if the vertices of $A$ and $B$ can be numbered so
that
$$a_0\le b_0\le a_1\le b_1\le a_2\le b_2\le a_3\le b_3\le a_0$$
where $a_i$, $0\le i\le 3$, are vertices of $A$ and $b_i$,
$0\le i\le 3$ are vertices of $B$.
\end{dfn}

Strong linkage is a closed condition: if two variable critical \ql s
are strongly linked and converge, then they must converge to two
strongly linked critical quadrilaterals. An obvious case of strong linkage is
between two non-degenerate critical quadrilaterals, whose vertices alternate
on the circle so that all the inequalities in
Definition~\ref{d:strolin} are strict. Yet even if both critical \ql
s are non-degenerate, some inequalities may be non-strict which
means that some vertices of both \ql s may coincide. For example,
two coinciding critical leaves can be viewed as strongly linked
critical \ql s, or an all-critical triangle $A$ with vertices $x, y,
z$ and its edge $B=\ol{yz}$ can be viewed as strongly linked \ql s
if the vertices are chosen as follows: $a_0=x, a_1=a_2=y, a_3=z$ and
$b_0=b_1=y, b_2=b_3=z$. If a critical \ql{} $Q$ is a critical leaf
or has all vertices distinct, then $Q$ as a critical \ql{} has a
well-defined set of vertices; the only ambiguous case is when $Q$ is
an all-critical triangle.

To study collections of critical \ql s we need a few notions and a
lemma. If a few chords can be concatenated to form a Jordan curve,
or if there are two identical chords, then we say that they form a
\emph{loop}. In particular, one chord does not form a loop while two
equal chords do. If an ordered collection of chords $(\ell_1, \dots,
\ell_k)$ contains no chords forming a loop we call it a \emph{no
loop collection}.

\begin{lem}\label{l:noloop}
The family of no loop collections of critical chords is closed.
\end{lem}

\begin{proof}
Suppose that there is a sequence of no loop collections of critical
chords $\nc^i=(\ell_1^i, \dots, \ell_s^i)$ with $\nc^i\to
\nc=(\ell_1, \dots, \ell_s)$ where all chords $\ell_i$ are critical.
We need to show that $\nc$ is a no loop collection. By way of
contradiction assume that, say, chords $\ell_1=\ol{a_1a_2}, \dots,
\ell_k=\ol{a_ka_1}$ form a loop $\widehat \nc$ in which the order of points $a_1,
\dots, a_k$ is positive. We claim that $\widehat \nc$ cannot be the limit
of no loop collections of critical chords, contradicting
the convergence assumption that $\nc^i\to \nc$.
This follows from the fact that if $G'\subset\uc$ is a union of
finitely many sufficiently small circle arcs such that all edges of
the convex hull $G=\ch(G')$ are critical, then in fact all circle
arcs in $G'$ are degenerate, so that $G$ is a finite polygon.
\end{proof}

Call a no loop collection of $d-1$ critical chords a \emph{full
collection}. Given a collection $\mathcal Q$ of $d-1$ critical \ql s
of a geolamination $\lam$, we choose one spike in each of them and
call this collection of $d-1$ critical chords a \emph{complete
sample of spikes (of $\mathcal Q$)}. If $\lam$ corresponds to a
lamination whose critical sets are critical quadrilaterals, any
complete sample of spikes is a full collection because in this case
distinct critical sets are disjoint. The fact that complete samples
of spikes form a full collection survives limit transition (unlike
pairwise disjointness). This inspires another definition.

\begin{dfn}[Quadratic criticality]\label{d:quacintro}
Let $(\lam, \qcp)$ be a geolamination with a $(d-1)$-tuple $\qcp$ of
critical \ql s that are gaps or leaves of $\lam$ such that any
complete sample of spikes is a full collection. Then $\qcp$ is
called a \emph{quadratically critical portrait $($qc-portrait$)$ for
$\lam$}
while the pair $(\lam, \qcp)$ is
called a \emph{geolamination with qc-portrait} (if the appropriate
geolamination $\lam$ for $\qcp$ exists but is not emphasized we
simply call $\qcp$ a \emph{qc-portrait}). The space of all
qc-portraits is denoted by $\fqcp_d$.
The family of all geolaminations with qc-portraits is denoted by
$\L\fqcp_d$.
\end{dfn}

If $C$ is
a complementary component of a complete sample of spikes in $\disk$, then $\si_d$
is one-to-one on the boundary of $C$ except for critical
chords contained in the boundary of $C$.

\begin{cor}\label{c:qcpoclose}
The spaces $\fqcp_d$
and $\L\fqcp_d$ are compact.
\end{cor}

\begin{proof}
Let $(\lam^i, \qcp^i)\to (\lam, \cp)$; by Theorem~\ref{t:sibliclos} and Lemma~\ref{l:qls}
here in the limit we have an invariant geolamination $\lam$ and an
ordered collection $\cp$ of $d-1$ critical \ql s. Let $\cp=(C_j)^{d-1}_{j=1}$ be the limit
critical \ql s. Choose a collection of spikes $\ell_j$ of \ql s of
$\cp$. Suppose that there is a loop formed by some of these spikes.
By construction there exist collections of spikes from
qc-portraits $\qcp^i$ converging to $(\ell_1, \dots, \ell_{d-1})$.
Since by definition these are full collections of critical chords,
this contradicts Lemma~\ref{l:noloop}. Hence $(\ell_1, \dots,
\ell_{d-1})$ is a full collection of critical chords too which
implies that $\cp$ is a qc-portrait for $\lam$ and proves that
$\fqcp_d$ and $\L\lp_d$ are compact spaces.
\end{proof}

The following lemma describes geolaminations admitting a
qc-portrait. Recall that by a \emph{collapsing \ql{}} we mean a critical
\ql{} which maps to a non-degenerate leaf.

\begin{lem}\label{l:qcpoexi}
A geolamination $\lam$ has a qc-portrait if and only if all
its critical sets are collapsing \ql s or all-critical sets.
\end{lem}

\begin{proof}
If $\lam$ has a qc-portrait, then the claim of the lemma follows by
definition.
Assume that the critical sets of $\lam$ are collapsing \ql s and
all-critical sets. Then $\lam$ may have several critical leaves.
Choose a maximal by cardinality no loop collection of critical
leaves of $\lam$. Add to them the collapsing \ql s of $\lam$.
Include all selected sets in the family of pairwise distinct sets
$\cp=(C_1, \dots, C_m)$ consisting of critical leaves and collapsing
\ql s. 


We claim that $\cp$ is a qc-portrait. To this end we need to show
that $m=d-1$ and that any collection $\nc$ of spikes of sets from
$\cp$ is a no loop collection. First let us show that any such
collection $\nc$ contains no loops. Indeed, suppose that $\nc$
contains a loop $\ell_1\in C_1$, $\dots$, $\ell_r\in C_r$. By
construction there must be a collapsing \ql{} among sets $C_1, \dots, C_r$.
We may assume that, say, $C_1=[a, x, b, y]$ is a collapsing
\ql{} 
and $\ell_1=\ol{ab}$ is
contained in the interior of $C_1$ except for points $a$ and $b$.
The spikes $\ell_2, \dots, \ell_r$ form a chain of concatenated
critical chords which has, say, $b$ as its initial point and $a$ as
its terminal point. Since these spikes come from sets $C_2, \dots,
C_r$ distinct from $C_1$, they have to pass through either $x$ or
$y$ as a vertex, a contradiction with $C_1$ being collapsing. 
Thus, $\nc$ contains no loops
which implies that the number $m$ of chords in $\nc$ is at most
$d-1$.

Assume now that $m<d-1$ and bring it to a contradiction. Indeed, if
$m<d-1$ then we can find a component $U$ of $\disk\sm \nc^+$ with
boundary including some circle arcs such that $\si_d$ on the
boundary of $\U$ is $k$-to-$1$ or higher with $k>1$ (images of critical
edges of $U$ may have more than $k$ preimages). We claim that there exists a critical
chord $\ell$ of $\lam$ inside $U$ that connects points in
$\bd(U)$ not connected by a chain of critical edges in
$\bd(U)$.
Observe that an arc on $\bd(U)$ may include several critical chords
from $\nc$. Consider all arcs $A\subset \bd(U)$ such that $\si_d$ is
strictly non-monotone on $A$, and the endpoints of $A$ are connected by a
leaf of $\lam$. Call such arcs \emph{non-monotone}.
Non-monotone arcs exist: by the assumptions there
exist leaves $\ell$ of $\lam$ inside $U$, and at least one of the two arcs in
the boundary of $U$ which connects the endpoints of $\ell$ must be
non-monotone.

The intersection of a decreasing sequence of non-monotone arcs is a
closed arc $A_0$ with endpoints connected with a leaf $\ell_0\in
\lam$ such that either $\ell_0$ is the desired critical leaf of
$\lam$ ($\ell_0$ cannot connect two points otherwise connected by a
chain of critical edges from $\bd(U)$ as this would contradict the
fact that arcs approaching $A_0$ are non-monotone), or $A_0$ is
still non-monotone. Thus, it is enough to show that if $A_0$ is a
minimal by inclusion non-monotone arc $A_0$ then there exists the
desired critical chord of $\lam$. 

Clearly, $A_0\cup \ell_0$ is a Jordan curve enclosing a Jordan disk
$T$, and $A_0$ is not a union of spikes. If $\ell_0$ is not critical
then by the assumption of minimality of $A_0$ the leaf $\ell_0$
cannot be approached by leaves of $\lam$ from within $T$, thus
$\ell_0$ is an edge of a gap $G\subset T$. Take a component $W$ of
$T\sm G$ which shares an edge $\m$ with $G$. Then, by minimality of
$A_0$, either $\bd(W)$ collapses to a point or $\bd(W)$ maps in a
non-strictly monotone fashion to the hole of $\si_d(G)$ located
``behind'' $\si_d(\m)$ united with $\si_d(\m)$. This implies that
$G$ is critical as otherwise the quoted properties of components $W$
of $T\sm G$ and the fact that $\si_d$ maps $G$ onto $\si_d(G)$ in a
one-to-one fashion show that $\si_d|_{A_0}$ is (non-strictly)
monotone, a contradiction. The gap $G$ cannot be all-critical, since
$\ell_0$ is an edge of $G$. Therefore, $G$ is a collapsing
quadrilateral, which contradicts our choice of $\Cc$.
\end{proof}

Observe that there might exist several qc-portraits for $\lam$ from
Lemma~\ref{l:qcpoexi}. For example, consider a $\si_4$-invariant
geolamination $\lam$ with two all-critical triangles $T_1=\ch(a,b,c),
T_2=\ch(a,c, d)$ sharing an edge $\ell=\ol{ac}$. The proof of
Lemma~\ref{l:qcpoexi} leads to a qc-portrait consisting of any three
edges of $T_1$, $T_2$ not equal to $\ell$ in some order (recall that for
each critical leaf its structure as a \ql{} is unique). However it is
easy to check that the collection $([a,b,b,c]$, $[a,a,c,c]$, $[a,c,d,
d])$ is a qc-portrait too.
Notice that, in the definition of a complete sample of spikes, we do not
allow to use more than one spike from each critical set, hence the fact
that the same spike appears twice in  $[a,a,c,c]$ does not result into
a loop.

Given a qc-portrait $\qcp$, any complete sample of spikes is a full
collection of critical chords. If $\qcp$ includes sets which are not
leaves, there are several complete samples of spikes as the choice
of spikes is ambiguous. This is important for Subsection~\ref{ss:smart},
where we introduce and study the so-called \emph{smart criticality}
and its applications to \emph{linked geolaminations with
qc-portraits} introduced below. First we need a technical
definition.

\begin{dfn}\label{d:compat}
A \emph{critical cluster} of $\lam$ is a
maximal by inclusion
convex subset of $\ol\disk$, whose boundary is a union of
critical leaves of $\lam$.
\end{dfn}

Consider the example
discussed after Lemma~\ref{l:qcpoexi}. There, a $\si_4$-invariant
geolamination $\lam$ has two all-critical triangles sharing a
critical edge; the union of these triangles is a
critical cluster
of $\lam$.

\begin{dfn}[Linked geolaminations]\label{d:qclink1}
Let $\lam_1$ and $\lam_2$ be geola\-mi\-na\-tions with qc-portraits
$\qcp_1=(C^i_1)_{i=1}^{d-1}$ and $\qcp_2=(C^i_2)_{i=1}^{d-1}$ and a
number $0\le k\le d-1$ such that:

\begin{enumerate}
\item for each $j>k$ the sets $C^j_1$ and $C^j_2$ are contained in
    a common critical cluster of $\lam_1$ and $\lam_2$ (in what
    follows these clusters will be called \emph{special critical clusters}
    and leaves contained in them will be called \emph{special
    critical leaves}).
\item for every $i$ with $1\le i\le k$, the sets $C^i_1$ and
    $C^i_2$ are either strongly linked critical quadrilaterals or share a
    spike.
\end{enumerate}

\noindent Then we use the following terminology:

\begin{enumerate}

\item[(a)] if in {\rm (1)} for every $i$ with $1\le i\le k$,
    the quadrilaterals $C^i_1$ and $C^i_2$ share a spike, we say that $\qcp_1$ and $\qcp_2$,
    (as well as $(\lam_1,
    \qcp_1)$ and $(\lam_2, \qcp_2)$) \emph{coincide in essence} (or
    \emph{essentially coincide}, or are \emph{essentially equal}),

\item[(b)] if in {\rm (1)} there exists $i$ with $1\le i\le k$ such
    that the \ql s $C^i_1$ and $C^i_2$ are strongly linked and do not share a
    spike, we say that $\qcp_1$ and $\qcp_2$ (as well as
    $(\lam_1, \qcp_1)$ and $(\lam_2, \qcp_2)$) are \emph{linked}.

\end{enumerate}

\noindent The critical sets $C^i_1$ and $C_2^i$, $1\le i\le d-1$ are
called \emph{associated $($critical sets of geolaminations with
qc-portraits $(\lam_1, \qcp_1)$ and $(\lam_2, \qcp_2))$}.
\end{dfn}

\subsection{Some special types of geolaminations}\label{ss:pandp}
Below, we discuss
perfect geolaminations and dendritic geolaminations.

\subsubsection{Perfect geolaminations}\label{sss:pergeo}

A geolamination $\lam$ is \emph{perfect} if no leaf of $\lam$ is
isolated. Every geolamination contains a maximal perfect
sublamination (clearly, this sublamination contains all degenerate
leaves). Indeed, consider $\lam$ as a metric space of leaves with the
Hausdorff metric and denote it by $\lam^*$. Then $\lam^*$ is a
compact metric space with a maximal perfect subset $\lam^c$
called \emph{the perfect sublamination of $\lam$}. The process of
finding $\lam^c$ was described in detail in \cite{bopt10}.
Lemma~\ref{l:perfect-sub-lam} follows from this description.

\begin{lem}\label{l:perfect-sub-lam}
The collection $\lam^c$ is an invariant perfect geolamination. For every $\ell\in
\lam^c$ and every neighborhood $U$ of $\ell$, there exist uncountably many
leaves of $\lam^c$ in $U$.
\end{lem}

Observe that there are at most two leaves of $\lam^c$ coming out of one point.
Otherwise, since, by Lemma~\ref{l:cones}, there are at most
countably many leaves of $\lam^c$ sharing an endpoint, $\lam^c$ has isolated leaves,
a contradiction.
Therefore,
any leaf of $\lam^c$ is a limit of an uncountably many leaves of
$\lam^c$ disjoint from $\ell$. If $\ell$ is critical, this implies
that $\si_d(\ell)$ is a point separated from the rest of the circle
by images of those leaves. Thus, a critical leaf $\ell$ is either
disjoint from all other leaves or gaps of $\lam^c$ or is an edge of
an all-critical gap of $\lam^c$ disjoint from all other leaves or
gaps of $\lam^c$. Together with the fact that at most
two leaves come out of a point, this implies Lemma~\ref{l:cridisj}.

\begin{lem}\label{l:cridisj}
Let $\lam$ be a perfect geolamination. Then the critical sets of
$\lam$ are pairwise disjoint and are either
all-critical sets, or critical sets mapping exactly $k$-to-$1$,
$k>1$, onto their images.
\end{lem}

\subsubsection{Dendritic geolaminations with critical patterns}\label{ss:fulden}
The main applications of our results will concern \emph{dendritic
laminations} defined below.

\begin{dfn}\label{d:dendrilam}
A lamination $\sim$ and its geolamination $\lam_\sim$ are called
\emph{dendritic} if the topological Julia set $J_\sim$ is a
dendrite. The family of all dendritic geolaminations is denoted by
$\L\mathcal D_d$.
\end{dfn}

Lemma~\ref{l:perfectd} is well-known.

\begin{lem}\label{l:perfectd}
Dendritic geolaminations $\lam$ are perfect.
\end{lem}

Dendritic geolaminations are closely related to polynomials. Let
$\dc$ be the space of all polynomials with connected Julia sets and
only repelling periodic points, and $\dc_d$ be the space of all such polynomials of degree
$d$. By Jan Kiwi's results \cite{kiwi97}, if a polynomial $P$ with
connected Julia set $J(P)$ has no Siegel or Cremer periodic points
(i.e., \emph{irrationally indifferent} periodic points whose
multiplier is of the form $e^{2\pi i\theta}$ for some irrational
$\theta$), then there exists a special lamination $\sim_P$,
determined by $P$, with the following property: $P|_{J(P)}$ is
monotonically semiconjugate to $f_{\sim_P}|_{J_{\sim_P}}$. Moreover,
all $\sim_P$-classes are finite, and the semiconjugacy is one-to-one
on all (pre)periodic points of $P$. These results apply to polynomials from $\dc$.

Strong conclusions about the topology of the Julia sets of
non-re\-nor\-ma\-li\-zable polynomials $P\in \dc$ follow from
\cite{kvs06}. Building upon earlier results by Kahn and Lyubich
\cite{kl09a, kl09b} and by Kozlovskii, Shen and van Strien
\cite{ksvs07a, ksvs07b}, Kozlovskii and van Strien generalized
results of Avila, Kahn, Lyubich and Shen \cite{akls09} and proved in
\cite{kvs06} that if all periodic points of $P$ are repelling, and
$P$ is non-renormalizable, then $J(P)$ is locally connected;
moreover, by \cite{kvs06}, two such polynomials that are
topologically conjugate are in fact quasi-conformally conjugate.
Thus, in this case $f_{\sim_P}|_{J_{\sim_P}}$ is a precise model of
$P|_{J(P)}$. Finally, for a given dendritic lamination $\sim$,
it follows from another result of Jan Kiwi \cite{kiw05} that there
exists a polynomial $P$ with $\sim=\sim_P$. Thus, by \cite{kiw05}
associating polynomials from $\dc$ with their laminations $\sim_P$
and geolaminations $\lam_P=\lam_{\sim_P}$, one maps polynomials from
$\dc_d$ \emph{onto} $\L\dc_d$.

To study the association of polynomials with their geolaminations, we
need Lemma~\ref{l:gm} (it is stated as a lemma in \cite{gm93} but
goes back to Douady and Hubbard \cite{hubbdoua85}).

\begin{lem}[\cite{gm93, hubbdoua85}]\label{l:gm}
Let $P$ be a polynomial, $\mathcal I$ be the set of all
$($pre$)$periodic external rays landing at the $P^n$-th preimage
$x_{-n}$ of a repelling periodic point $x$ so that $x_{-n}$ be
not $($pre$)$critical. Then the set $\mathcal I$ is finite, and for
any polynomial $P^*$ sufficiently close to $P$, there is a
corresponding repelling periodic point $x^*$ close to $x$ and there is a
$(P^*)^n$-th preimage $x^*_{-n}$ of $x^*$ close to $x_{-n}$ such
that the family $\mathcal I^*$ of all $($pre$)$periodic rays,
landing at $x^*_{-n}$, consists of rays uniformly $($with respect to
the spherical metric$)$ close to the
corresponding rays of $\mathcal I$ with the same external arguments.
\end{lem}

We also need the following 
lemma.

\begin{lem}\label{l:condense}
Suppose that $\sim$ is a dendritic lamination. Then each leaf of
$\lam_\sim$ can be approximated by $($pre$)$periodic leaves.
\end{lem}

\begin{proof}
Consider the topological polynomial $f_\sim$. Choose an arc
$I\subset J_\sim$. By \cite{bl02}, we can find $k>0$ such that $I$
and $f^k_\sim(I)$ are non-disjoint. Consider the union $T$ of all
$f^k_\sim$-images of $I$ (this union is connected) and take its closure
$K$. Then $K\subset J_\sim$ is an $f^k_\sim$-invariant dendrite. Any
periodic point $x\in K$ corresponds to a $\sim$-class whose convex hull has
periodic edges fixed by $\si_d^m$ for some $m>0$. Hence
there are short open pairwise disjoint arcs $(x, s')\subset (x,
s)\subset K$ such that all points $y\in (x, s')$ are repelled away
from $x$ but have images in $(x, s)$. By Theorem 7.2.6
of \cite{bfmot10}, there are infinitely many periodic cutpoints in
$K$. Since $T$ is connected and dense in $K$, it follows that $T$
contains periodic points. Hence $I$ contains (pre)periodic points.
Clearly, this implies the lemma.
\end{proof}

We will use qc-portraits to parameterize (tag) dendritic
geolaminations. An obstacle to this is the fact that a geolamination
$\lam$ with a $k$-to-$1$ critical set such that $k>2$ does not admit a
qc-portrait. However, using Lemma~\ref{l:qcpoexi}, it is easy
to see that in this case one can insert critical \ql s in critical sets
of higher degree in order to ``tune'' $\lam$ into a geolamination with
a qc-portrait. This motivates the following.

\begin{dfn}\label{d:cripatt1}
Let $\lam$ have pairwise disjoint critical sets (gaps or leaves) $D_1$, $\dots$, $D_k$. Let
$\lam\subset \lam_1$ and $\qcp=(E_1,\dots,E_{d-1})$ be a qc-portrait
for $\lam_1$. Clearly, there is a unique $(d-1)$-tuple $\zc$ $=$
$(C_1,\dots,C_{d-1})$ such that for every $1\le i\le d-1$ we have
$E_i\subset C_i$ and there is $1\le j(i)\le k$ with $C_i=D_{j(i)}$. Then $\zc$
is called the \emph{critical pattern of $\qcp$ in $\lam$}. Observe that
each $D_{j(i)}$ is repeated in $\zc$ exactly $m_{j(i)}-1$ times, where $m_{j(i)}$ is
the degree of $D_{j(i)}$.

In general, given a geolamination $\lam$ with pairwise disjoint
critical sets  $D_1, \dots, D_k$, by a
\emph{geolamination with a critical pattern} we mean a pair $(\lam,
\zc)$ where $\zc$ $=$ $(C_1,\dots,C_{d-1})$ is a $(d-1)$-tuple of sets
provided for every $1\le i\le d-1$ there is a $1\le j\le k$ with
$C_i=D_j$ and, for every $j=1$, $\dots$, $k$, each $D_j$ is repeated in
$\zc$ exactly $m_j-1$ times, where $m_j$ is the degree of $D_j$. Then
$\zc$ is called a \emph{critical pattern for $\lam$}. The space of all
\emph{dendritic} geolaminations with critical patterns is denoted by
$\L\cpd_d$.
\end{dfn}

By changing the order of the critical sets, various critical patterns
for the same geolamination can be obtained. In the dendritic case,
the connection between critical patterns and geolaminations can be
studied using results of Jan Kiwi \cite{kiwi97}. One of the results
of \cite{kiwi97} can be stated as follows: if $\lam$ is a dendritic
geolamination and $\lam'$ is an invariant geolamination such that
$\lam$ and $\lam'$ share a collection of $d-1$ critical chords with
no loops among them, then $\lam'\supset \lam$. Since all gaps of
$\lam$ are finite, this means that $\lam'\sm \lam$ consists of
countably many leaves inserted in certain gaps of $\lam$.

Observe that if a sequence of geolaminations with critical patterns
$(\lam^i, \zc^i)$ converges, then, by
Theorem~\ref{t:sibliclos}, the limit $\lam^\infty$ of geolaminations
$\lam^i$ is itself a $\si_d$-invariant geolamination. Moreover, it
is easy to see that then critical patterns $\zc^i$ converge to the
limit collection of $d-1$ critical sets of $\lam^\infty$. Together
with results from \cite{kiwi97}, this implies the following lemma.

\begin{lem}\label{l:uppers}
Suppose that a sequence of geolaminations with critical patterns $(\lam^i,
\zc^i)$ converges in the sense of the Hausdorff metric to a
geolamination $\lam^\infty$ with a collection of limit critical sets
$C_1, \dots, C_{d-1}$. Suppose that there exists a dendritic
geolamination $\lam$ with a critical pattern $\zc=(Z_1, \dots,
Z_{d-1})$ such that $C_i\subset Z_i, 1\le i\le d-1$. Then
$\lam^\infty\supset \lam$.
\end{lem}

For an integer $m>0$, we use a partial order by inclusion among
$m$-tuples: $(A_1, \dots, A_m)\succ (B_1, \dots, B_m)$ (or
$(B_1, \dots, B_m)\prec (A_1, \dots, A_m)$) if and only if
$A_i\supset B_i$ for all $i=1$, $\dots$, $m$. Thus $m$-tuples and
$k$-tuples with $m\ne k$ are always incomparable.
Lemma~\ref{l:uppers} says that if critical patterns converge
\emph{into} a critical pattern of a dendritic geolamination $\lam$,
then the corresponding geolaminations themselves converge
\emph{over} $\lam$.

The notion of a geolamination with critical pattern is related to
the notion of a \emph{$($critically$)$ marked polynomial} \cite{mil12},
i.e., a polynomial $P$ with an ordered collection $\cm$ of its
critical points, each of which is listed according to its
multiplicity (so that there are $d-1$ points in $\cm$). Critically
marked polynomials do not have to be dendritic (in fact, the notion
is used by Milnor and Poirier for hyperbolic polynomials, i.e., in
the situation diametrically opposite to that of dendritic
polynomials). Evidently, the space of critically marked polynomials
is closed, and if $P$ is perturbed a little, the critical points of
the perturbed polynomial can be ordered to give rise to a critically
marked polynomial close to the original $(P, \cm)$ (that is, the
natural forgetful map from critically marked polynomials to
polynomials is a branched covering).

Denote the space of all degree $d$ critically marked dendritic
polynomials by $\cmd_d$. To each $(P, \cm)\in \cmd_d$ we associate
the corresponding dendritic geolamination with a critical pattern
$(\lam_{\sim_P}, \zc)$ in a natural way (each point $z\in J(P)$ is
by \cite{kiwi97} associated to a gap or leaf $G_z$ of
$\lam_{\sim_P}$, thus each point $c\in \cm$ is associated with the
critical gap or leaf $G_c$ of $\lam_{\sim_P}$). This defines the map
$\Psi_d: \cmd_d\to \L\cpd_d$ such that $\Psi_d(P, \cm)=(\lam_{\sim_P},
\zc)$. Corollary~\ref{c:crista} easily follows
from Lemmas~~\ref{l:gm}, ~\ref{l:condense} and ~\ref{l:uppers}.

\begin{cor}\label{c:crista}
Suppose that a sequence $(P_i, \cm_i)$ of critically marked dendritic
polynomials converges to a critically marked dendritic polynomial $(P, \cm)$.
Set $(\lam_{\sim_{P_i}}, \zc_i)=\Psi_d(P, \cm_i)$ and
$(\lam_{\sim_P}, \zc)=\Psi_d(P, \cm)$.
If $(\lam_{\sim_{P_i}}, \zc_i)$ converge in the sense of the Hausdorff
metric to $(\lam^\infty,\zc_\infty)$, then
$\lam^\infty\supset \lam_{\sim_P}$ and $\zc_\infty\prec \zc$.
\end{cor}

By Corollary~\ref{c:crista}, critical sets of geolaminations
$\lam_{\sim_P}$ associated with polynomials $P\in \dc_d$ cannot explode
under perturbation of $P$ (they may implode though). Provided a
geometric (visual) way to parameterize $\L\cpd_d$ is given, the map
$\Psi_d$ yields the corresponding parameterization of $\cmd_d$ and
gives an important application of our tools.

\section{Accordions of Laminations}\label{s:acclam}

In the Introduction, we mentioned that some of Thurston's tools from
\cite{thu85} fail in the cubic case. This motivates us to develop new
tools (so-called \emph{accordions}), which basically
track linked leaves from different geolaminations. In this
section, we study accordions in detail. In Sections~\ref{s:acclam} and
\ref{s:qcrit}, we assume that $\lam_1$, $\lam_2$ are $\si_d$-invariant
geolaminations, and $\ell_1$, $\ell_2$ are leaves of $\lam_1$, $\lam_2$,
respectively.

\subsection{Motivation}\label{ss:motiv}
For a quadratic invariant geolamination $\lam$ and a leaf $\ell$ of $\lam$
that is not a diameter, let $\ell'$ be the sibling of $\ell$
(disjoint from $\ell$). Denote by $C(\ell)$ the open strip of
$\cdisk$ between $\ell$ and $\ell'$ and by $L(\ell)$ the length of
the shorter component of $\uc\sm \ell$.
Suppose that $\frac13\le L(\ell)<\frac12$, and that $k$ is the
smallest number such that $\si_2^k(\ell)\subset C(\ell)$ except
perhaps for the endpoints. The Central Strip Lemma (Lemma II.5.1 of
\cite{thu85}) claims that $\si_2^k(\ell)$ separates $\ell$ and
$\ell'$. In particular, if $\ell=M$ is a \emph{major}, i.e., a
longest leaf of some quadratic invariant geolamination, then an eventual
image of $M$ cannot enter $C(M)$.

Let us list Thurston's results for which the Central Strip Lemma is
crucial. A \emph{$\si_2$-wandering triangle} is a triangle with
vertices $a$, $b$, $c$ on $\uc$ such that the convex hull
$T_n$ of $\si_2^n(a)$, $\si_2^n(b)$, $\si_2^n(c)$ is a
non-degenerate triangle for every $n=0$, $1$, $\dots$, and all these
triangles are pairwise disjoint.

\begin{thm}[No Wandering Triangle Theorem \cite{thu85}]\label{t:nwt}
There are no wandering triangles for $\si_2$.
\end{thm}

Theorem~\ref{t:gaptrans} stated below follows from the Central Strip
Lemma and is due to Thurston for $d=2$. For arbitrary $d$, it is due
to Jan Kiwi, who used different tools.

\begin{thm}[\cite{thu85, kiw02}]\label{t:gaptrans}
If $A$ is a finite $\si_d$-periodic gap of period $k$, then either
$A$ is a $d$-gon, and $\si^k_d$ fixes all vertices of $A$, or there
are at most $d-1$ orbits of vertices of $A$ under $\si_d^k$. Thus, for
$d=2$, the remap is transitive on the vertices of any finite
periodic gap.
\end{thm}

\begin{figure}
  \includegraphics[width=4.5cm]{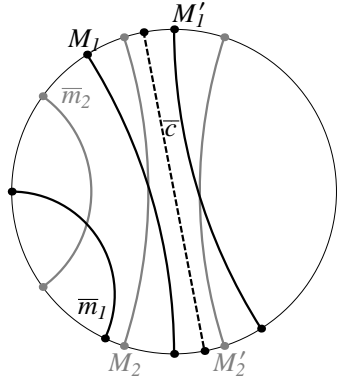}
  \caption{This figure illustrates Thurston's proof that quadratic minors are unlinked.
  The Central Strip Lemma forces orbits of both minors to not cross $\ovc$.}
  \label{f:central-strip}
\end{figure}

Another crucial result of Thurston is that minors of distinct
quadratic invariant geolaminations are disjoint in $\disk$. A sketch
of the argument follows. Let $\m_1$ and $\m_2$ be the minors of two
invariant geolaminations $\lam_1\ne \lam_2$ that cross
in $\disk$. Let $M_1$, $M'_1$ and $M_2$, $M'_2$ be the two pairs of
corresponding majors. We may assume that $M_1$, $M_2$ cross in $\disk$ and
$M'_1$, $M'_2$ cross in $\disk$,
but $(M_1\cup M_2)\cap (M'_1\cup M'_2)=\0$ (see Figure
\ref{f:central-strip}) so that there is a diameter $\ovc$ with
strictly preperiodic endpoints separating $M_1\cup M_2$ from
$M_1'\cup M_2'$. Thurston shows that there is a unique invariant
geolamination $\lam$, with only finite gaps, whose major is $\ovc$.
By the Central Strip Lemma, forward images of $\m_1$, $\m_2$ do not
intersect $\ovc$. Hence $\m_1\cup \m_2$ is contained in
a finite gap $G$ of $\lam$. By the No Wandering Triangle Theorem, $G$
is eventually periodic. By Theorem~\ref{t:gaptrans}, some images of
$\m_1$ intersect inside $\disk$, a contradiction.

\begin{figure}
\includegraphics[width=4.5cm]{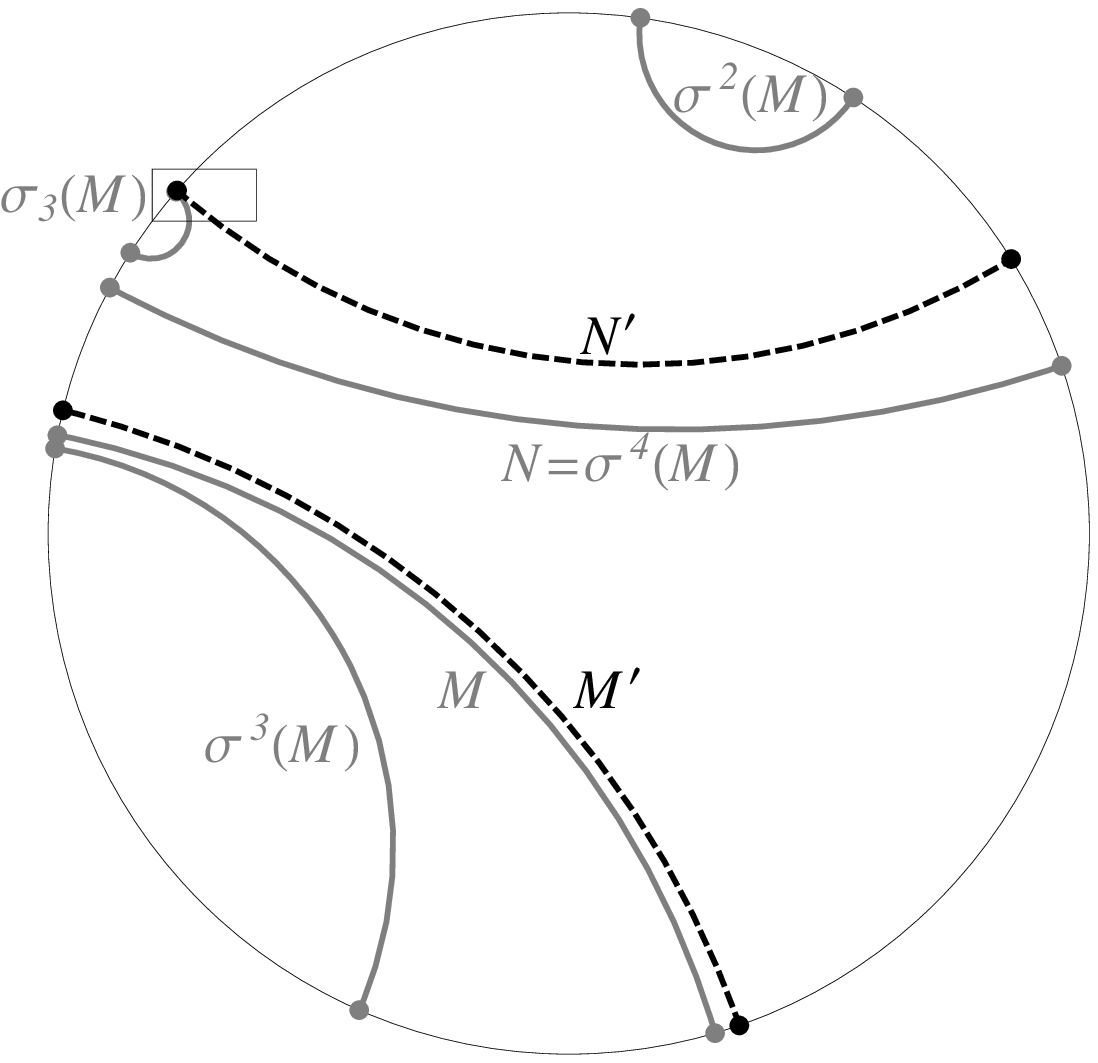}
\hspace{2cm}
\includegraphics[width=4.5cm]{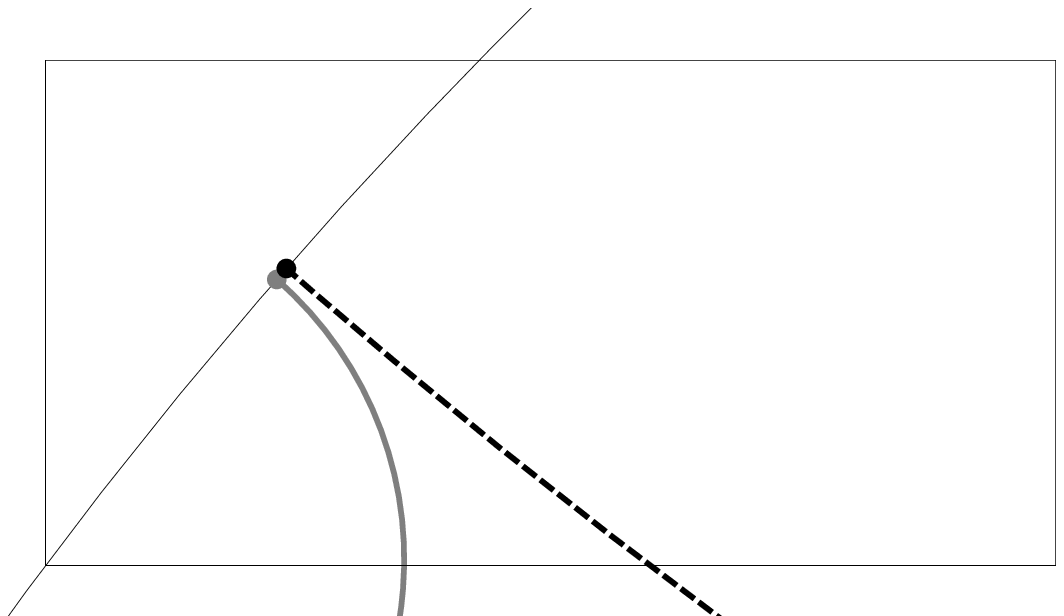}
{\caption{This figure shows that the Central Strip Lemma fails in the
cubic case. Its left part has a fragment in which two endpoints of
leaves are located very close to each other. Its right part is
the zoomed-in version of the fragment indicating that the periodic
points do not coincide.}\label{f:ex-intro}}
\end{figure}

Examples indicate that statements analogous to the Central Strip
Lemma fail in the cubic case. Indeed,
Figure~\ref{f:ex-intro} shows a leaf
$M=\ol{\frac{342}{728}\frac{579}{728}}$ of period $6$ under $\si_3$
and its $\si_3$-orbit together with the leaf $M'$ (which has the
same image as $M$ forming together with $M$ a narrower ``critical
strip'' $S_n$) and the leaf $N'$ (which has the same image as
$N=(\si_3)^4(M)$ forming together with $N$ a wider ``critical
strip'' $S_w$). Observe that $\si_3(M)\subset S_w$, which shows that
the Central Strip Lemma does not hold in the cubic case (orbits of
periodic leaves may give rise to ``critical strips'' containing some
elements of these orbits of leaves). This apparently makes a direct
extension of the arguments from the previous paragraph impossible
leaving the issue of whether and how minors of cubic geolaminations
can be linked unresolved.

Another consequence of the failure of the Central Strip Lemma in the
cubic case is the failure of the No Wandering Triangle Theorem (a
counterexample was given in \cite{bo08}). Properties of wandering
polygons were studied in \cite{kiw02, bl02, chi07}.

\subsection{Properties of accordions}\label{ss:propacc}
We now give the definition of accordions.
\begin{dfn}\label{d:accord}
Let $A_{\lam_2}(\ell_1)$ be the collection of leaves of $\lam_2$
linked with $\ell_1$, together with $\ell_1$. Let
$A_{\ell_2}(\ell_1)$ be the collection of leaves from the forward
orbit of $\ell_2$ that are linked with $\ell_1$, together with
$\ell_1$. The sets defined above are called \emph{accordions $($of
$\ell_1)$}
while $\ell_1$ is called the \emph{axis} of the accordion.
Sometimes we will also use $A_{\lam_2}(\ell_1)$ and $A_{\ell_2}(\ell_1)$
to mean the union of the leaves constituting these accordions.
\end{dfn}

In general, accordions do not behave nicely
under $\si_d$
as leaves which are linked may have unlinked images. To avoid these
problems, for the rest of this section, we will impose the following
conditions on accordions.

\begin{dfn}\label{d:opacc}
A leaf $\ell_1$ is said to \emph{have order preserving accordions
with respect to $\lam_2$ $($respectively, to a leaf $\ell_2)$} if
$A_{\lam_2}(\ell_1)\ne \{\ell_1\}$ (respectively, $A_{\ell_2}(\ell_1)\ne
\{\ell_1\}$), and, for each $k\ge 0$, the map $\si_d$ restricted to
$A_{\lam_2}(\si_d^k(\ell_1))\cap\uc$ (respectively, to
$A_{\ell_2}(\si_d^k(\ell_1))\cap\uc$)
is order preserving (in particular, it is one-to-one). Say that $\ell_1$ and $\ell_2$ have
\emph{mutually order preserving accordions} if $\ell_1$ has order
preserving accordions with respect to $\ell_2$, and vice versa
(in particular, $\ell_1$ and $\ell_2$ are not precritical).
\end{dfn}

Though fairly strong, these conditions naturally arise in
the study of linked or essentially coinciding geolaminations.  In Section~\ref{s:qcrit}, we
will show that they are often satisfied by pairs of linked leaves
of linked or essentially coinciding geolaminations (Lemma~\ref{l:indepcrit}) so that
there are at most countably many pairs of linked leaves which do not have
mutually order preserving accordions. If
geolaminations are perfect, this will imply that every accordion
consisting of more than one leaf contains a pair of leaves with
mutually order preserving accordions. Understanding the rigid
dynamics of such pairs is crucial to our main results.

The proof of
Proposition~\ref{p:accimage} is left to the reader.

\begin{prop}\label{p:accimage}
If $\si_d$ is order preserving on an accordion $A$ with axis
$\ell_1$ and $\ell\in A$, $\ell\ne\ell_1$, then $\si_d(\ell)$ and
$\si_d(\ell_1)$ are linked. In particular, if $\ell_1$ has order
preserving accordions with respect to $\ell_2$ then
$\si_d^k(\ell)\in A_{\ell_2}(\si_d^k(\ell_1))$ for every $\ell\in
A_{\ell_2}(\ell_1)$, $\ell\ne\ell_1$, and every $k\ge 0$.
\end{prop}

We now explore more closely the orbits of leaves from
Definition~\ref{d:opacc}.

\begin{prop}\label{p:sameori}
Suppose that $\ell_1$ and $\ell_2$ are linked, $\ell_1$ has order
preserving accordions with respect to $\ell_2$, and $\si_d^k(\ell_2)\in
A_{\ell_2}(\ell_1)$ for some $k>0$. In this case, if $\ell_2=\ol{xy}$,
then either $\ell_1$ separates $x$ from $\si_d^k(x)$ and $y$ from
$\si_d^k(y)$, or $\ell_2$ has $\si_d^k$-fixed endpoints.
\end{prop}

\begin{figure}[H]
\includegraphics[width=4.5cm]{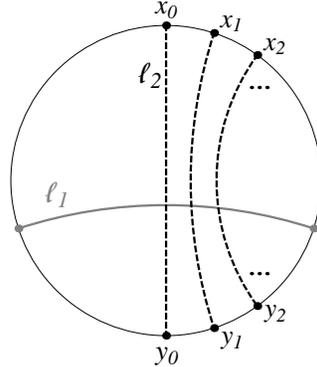}
{\caption{This figure illustrates
Proposition~\ref{p:sameori}. Although in the figure $\ol{x_2y_2}$
is linked with $\ell_1$, the argument does not assume this.
In this and forthcoming figures, leaves marked in the same fashion
belong to the same grand orbits of leaves.
}\label{f:sameori}}
\end{figure}

\begin{proof}
Suppose that $\ell_2$ is not $\si_d^k$-fixed. Denote by
$x_0=x,y_0=y$ the endpoints of $\ell_2$; set $x_i=\si_d^{ik}(x_0),
y_i=\si_d^{ik}(y_0)$ and $A_t=A_{\ell_2}(\si_d^t(\ell_1)), t=0, 1,
\dots.$ If $\ell_1$ does not separate $x_0$ and $x_1$, then either
$x_0\le x_1<y_1\le y_0$ or $x_0<y_0\le y_1<x_1\le x_0$. We may
assume the latter (cf. Figure \ref{f:sameori}).

Since $\si_d^k$ is order preserving on $A_0\cap\uc$, then $x_0<y_0\le
y_1\le y_2<x_2 \le x_1 \le x_0$ while the leaves $\ol{x_1y_1}$ and
$\ol{x_2y_2}$ belong to the accordion $A_k$ so that the above
inequalities can be iterated. Inductively we see that

$$x_0< y_0 \le \dots \le y_{m-1}\le y_m<x_m\le x_{m-1}\le \dots \le x_0.$$

\noindent All leaves $\ol{x_i y_i}$ are pairwise distinct as
otherwise there exists $n$ such that $\ol{x_{n-1} y_{n-1}}\ne
\ol{x_n y_n}=\ol{x_{n+1} y_{n+1}}$ contradicting $\si_d^k$ being
order preserving on $A_{k(n-1)}$. Hence the leaves $\ol{x_i y_i}$
converge to a $\si_d^k$-fixed point or leaf, contradicting
the expansion property of $\si_d^k$.
\end{proof}

In what follows, we often use one of the endpoints of a leaf
as the subscript in the notation for this leaf.

\begin{lem}\label{l:sameperiod}
If $\ell_a=\ol{ab}$ and $\ell_x=\ol{xy}$, where $a<x<b<y$, are linked
leaves with mutually order preserving accordions, and $a$, $b$ are of
period $k$, then $x$, $y$ are also of period $k$.
\end{lem}

\begin{proof}
By the order preservation, $\si_d^k(x)$ is not separated from $x$ by $\ell_a$.
It follows from Proposition~\ref{p:sameori} that $x=\si^k_d(x)$, $y=\si^k_d(y)$.
Since, by Lemma~\ref{l:sameperiod1}, the points $x$ and $y$ have the same
period (say, $m$), then $m$ divides $k$. Similarly, $k$ divides $m$.
Hence $k=m$.
\end{proof}

We will mostly use the following corollary of the above results.

\begin{figure}
\includegraphics[width=4.5cm]{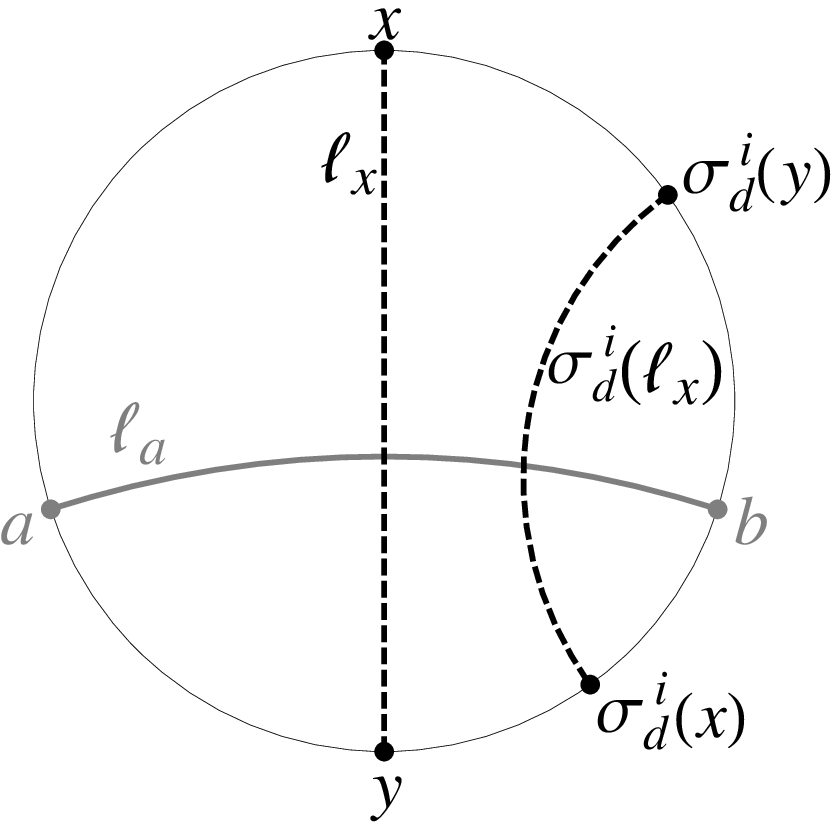}
\hspace{2cm}
\includegraphics[width=4.5cm]{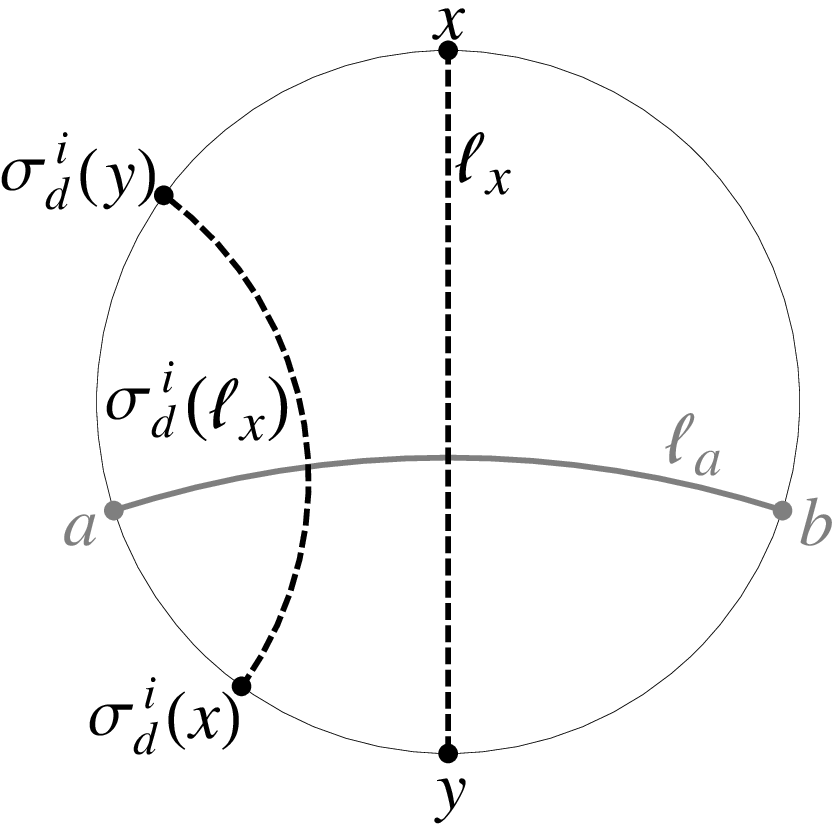}
{\caption{This figure shows two cases listed in
Corollary~\ref{c:nothreelink}, part (4).}\label{fig:3link}}
\end{figure}

\begin{cor}\label{c:nothreelink}
Suppose that $\ell_a=\ol{ab}$ and $\ell_x=\ol{xy}$ with $x<a<y<b$
are linked leaves. If $\ell_a$ and $\ell_x$ have mutually order
preserving accordions, then there are the following possibilities for
$A=A_{\ell_x}(\ell_a)$.

\begin{enumerate}

\item $A=\{\ell_a, \ell_x\}$ and no forward image of $\ell_x$
crosses $\ell_a$.

\item $A=\{\ell_a, \ell_x\}$, the points $a$, $b$, $x$, $y$ are
of period $2j$ for some $j$,
$\si^j(x)=y, \si^j(y)=x$, and either $\si^j_d(a)=b$, $\si^j_d(b)=a$,
or $\si_d^j(\ell_a)\ne \ell_a$, and $\ell_x$ separates
the points $a$, $\si_d^j(b)$ from the points $b$, $\si_d^j(a)$.

\item $A=\{\ell_a, \ell_x\}$, the points $a$, $b$, $x$, $y$ are of the
same period, $x$, $y$ have distinct orbits, and $a$, $b$
have distinct orbits.

\item There exists $i>0$ such that $A=\{\ell_a, \ell_x,
\si^i_d(\ell_x)\}$ and either $x<a<y\le \si_d^i(x)<b<\si_d^i(y)\le
x$ or $x\le \si^i_d(y)<a<\si^i_d(x)\le y<b$, as shown in Figure \ref{fig:3link}.
\end{enumerate}

\end{cor}

\begin{proof}
Three \emph{distinct} images of $\ell_x$ cannot cross $\ell_a$ as if
they do, then it is impossible for the separation required in
Proposition~\ref{p:sameori} to occur for all of the pairs of images
of $\ell_x$. Hence at most two images of $\ell_x$ cross $\ell_a$.

If two distinct leaves from the orbit of $\ell_x$
cross $\ell_a$, then, by Proposition~\ref{p:sameori} and the order
preservation, case (4) holds.
Thus we can assume that $A=\{\ell_a, \ell_x\}$.
If no forward image of $\ell_x$ is linked with $\ell_a$, then we have case (1).

In all remaining cases we have $\si_d^k(\ell_x)=\ell_x$ for some
$k>0$. By Lemma~\ref{l:sameperiod1}, points $x$ and $y$ are of the
same period. Suppose that $x$, $y$ belong to the same periodic
orbit. Choose the least $j$ such that $\si^j_d(x)=y$. Let us show
that then $\si^j(y)=x$. Indeed, assume that $\si^j(y)\ne x$. Since
by the assumption the only leaf from the forward orbit of $\ell_x$,
linked with $\ell_a$, is $\ell_x$, we may assume (for the sake of
definiteness) that $y<\si^j_d(y)\le b$. Then a finite concatenation
of further $\si^j_d$-images of $\ell_x$ will connect $y$ with $x$.
Again, since $A=\{\ell_a, \ell_x\}$, one of their endpoints will
coincide with $b$. Thus, $y<\si^j_d(y)\le b<\si_d^j(b)\le x$, see
Figure \ref{f:nothreelinka2i}. Let us now apply $\si_d^j$ to $A$; by the
order preservation $y<\si^j_d(a)<\si^j_d(y)\le b<\si^j_d(b)\le x<a$.
Hence, $\si^j_d(\ell_a)$ is linked with $\ell_a$, a contradiction.

\begin{figure}
\includegraphics[width=4.5cm]{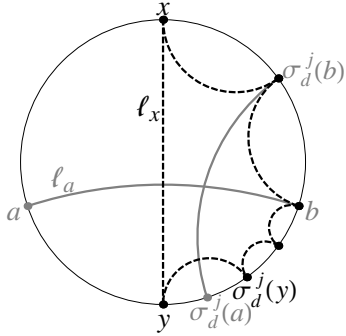}
{\caption{This figure illustrates the proof of
Corollary~\ref{c:nothreelink}.}\label{f:nothreelinka2i}}
\end{figure}

Thus, $\si^j(y)=x$ (i.e., $\si_d^j$ flips $\ell_x$ onto itself),
$k=j$, the points $x$ and $y$ are of period $2j$ and, by
Lemma~\ref{l:sameperiod}, the points $a$ and $b$ are also of period
$2j$. If $\si_d^j(a)=b$, then $\si_d^j(b)=a$, and
if $\si_d^j(b)=a$, then $\si_d^j(a)=b$ (since both points have period $2j$). Now, if
$\si_d^j(a)\ne b$ and $\si_d^j(b)\ne a$, then, by the order preservation,
$\ell_x$ separates the points $a$, $\si_d^j(b)$ from the points $b$, $\si_d^j(a)$.
So, case (2) holds.

Assume that $x$ and $y$ belong to distinct periodic
orbits of period $k$. By Lemma~\ref{l:sameperiod}, the points $a$, $b$ are
of period $k$. Let points $a$ and $b$ have the same orbit. Then, if
$k=2i$ and $\si_d^i$ flips $\ell_a$ onto itself, it would follow from the
order preservation that $\si_d^i(\ell_x)$ is linked with $\ell_a$.
Since $\ell_x$ is the unique leaf from the orbit of $\ell_x$ linked
with $\ell_a$ this would imply that $\si_d^i$ flips $\ell_x$ onto
itself, a contradiction with $x, y$ having disjoint orbits. Hence we
may assume that, for some $j$ and $m>2$, we have that $\si_d^j(a)=b$,
$jm=k$, and a concatenation of leaves $\ell_a$, $\si_d^j(\ell_a)$,
$\dots$, $\si_d^{j(m-1)}(\ell_a)$ forms a polygon $P$.

If one of these leaves distinct from $\ell_a$ (say,
$\si_d^{js}(\ell_a)$) is linked with $\ell_x$, we can apply the map
$\si_d^{j(m-s)}$ to $\si_d^{js}(\ell_a)$ and $\ell_x$; by order
preservation we will see then that $\ell_a$ and
$\si_d^{j(m-s)}(\ell_x)\ne \ell_x$ are linked, a contradiction with
the assumption that $A=\{\ell_a, \ell_x\}$. If none of the leaves
$\si_d^j(\ell_a)$, $\dots$, $\si_d^{j(m-1)}(\ell_a)$ is linked with
$\ell_x$, then $P$ has an endpoint of $\ell_x$ as one of its
vertices. As in the argument given above, we can then apply
$\si_d^j$ to $A$ and observe that, by the order preservation, the
$\si_d^j$-image of $\ell_x$ is forced to be linked with $\ell_x$, a
contradiction. Hence $a$ and $b$ have disjoint orbits, and case (3)
holds.
\end{proof}

\subsection{Accordions are (pre-)periodic or wandering}\label{ss:accwape}
Here we prove Theorem~\ref{t:compgap} which is the main result of
Section~\ref{s:acclam}.

\begin{dfn}\label{d:positor}
A finite sequence of points $x_0, \dots, x_{k-1}\in \uc$ is
\emph{positively ordered} if $x_0<x_1<\dots<x_{k-1}<x_0$. If the
inequality is reversed, then we say that points $x_0$, $\dots$,
$x_{k-1}\in \uc$ are \emph{negatively ordered}. A sequence $y_0$, $y_1$,
$\dots$ is said to be \emph{positively circularly ordered} if it is
either positively ordered or there exists $k$ such that $y_i=y_{i
\mod k}$ and $y_0<y_1<\dots<y_{k-1}<y_0$. Similarly we define points
that are \emph{negatively circularly ordered}.
\end{dfn}

A positively (negatively) \emph{circularly} ordered sequence that
is not positively (negatively) ordered is a sequence, whose points
repeat themselves after the initial collection of points that are
positively (negatively) ordered.

\begin{dfn}\label{d:leavineq}
Suppose that the chords $\bt_1,$ $\dots,$ $\bt_n$ are edges of the
closure $Q$ of a single component of $\disk\sm \bigcup \bt_i$. For each
$i$, let $m_i$ be the midpoint of the hole $H_Q(\bt_i)$. We write
$\bt_1<\bt_2<\dots<\bt_n$ if the points $m_i$ form a positively ordered
set and call the chords $\bt_1$, $\dots$, $\bt_n$ \emph{positively
ordered}. If the points $m_i$ are positively circularly
ordered, then we say that $\bt_1,$ $\dots,$ $\bt_n$ are
\emph{positively circularly ordered}. \emph{Negatively ordered} and
\emph{negatively circularly ordered} chords are defined similarly.
\end{dfn}

Lemma~\ref{l:linkstruct} is used in the main result of this section.

\begin{lem}\label{l:linkstruct} If $\ell_a$ and $\ell_x$
are linked, have mutually order preserving accordions, and
$\si_d^k(\ell_x)\in A_{\ell_x}(\ell_a)$ for some $k>0$, then, for
every $j>0$, the leaves $\si_d^{ki}(\ell_x)$, $i=0$, $\dots$, $j$, are
circularly ordered, and $\ell_a$, $\ell_x$ are periodic with
endpoints of the same period.
\end{lem}

\begin{figure}
\includegraphics[width=4.5cm]{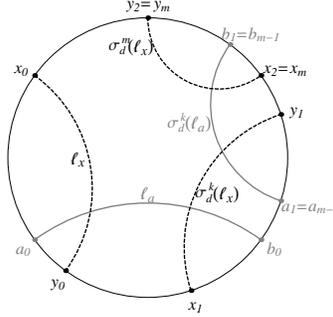}
{\caption{This figure illustrates Lemma~\ref{l:linkstruct}. Images
of $\ell_a$ cannot cross other images of $\ell_a$, neither can they
cross images of $\ell_x$ that are already linked with two images of
$\ell_a$ (by Corollary~\ref{c:nothreelink}). Similar claims hold for
$\ell_x$.}\label{f:linkstruct}}
\end{figure}

\begin{proof}
By Lemma~\ref{l:sameperiod}, we may assume that case (4) of
Corollary~\ref{c:nothreelink} holds (and so $\si_d^k(\ell_x)\ne
\ell_x$). Set $B=\{\ell_a, \ell_x\}$, $\ell_a=\ol{ab},
\ell_x=\ol{xy}$ and let $a_i$, $b_i$, $x_i$, $y_i$ denote the
$\si_d^{ik}$-images of $a$, $b$, $x$, $y$, respectively ($i\ge 0$). We may
assume that the first possibility from case (4) holds and
$x_0<a_0<y_0\le x_1<b_0<y_1\le x_0$ (see the left part of Figure
\ref{fig:3link} and Figure \ref{f:linkstruct}). By the assumption of mutually
order preserving accordions applied to $B$, we have $x_i<a_i<y_i\le
x_{i+1}<b_i<y_{i+1}\le x_i$ ($i\ge 0$), in particular $x_1<a_1<y_1$.
Then there are two cases depending on the location of $a_1$.
Consider one of them as the other one can be considered similarly.
Namely, assume that $b_0<a_1<y_1$ and proceed by
induction for $m$ steps observing that

$$x_0<a_0<y_0\le x_1<b_0\le a_1<\dots\le x_m<b_{m-1}<a_m<y_m\le x_0.$$

\noindent Thus, the first $m$ iterated $\si_d^k$-images of $\ell_x$ are
circularly ordered and alternately linked with the first $m-1$ iterated
images of $\ell_a$ under $\si_d^k$ (see Figure \ref{f:linkstruct}).
In the rest of the proof, we exploit the following fact.

\smallskip

\noindent\textbf{Claim A.} \emph{Further images of $\ell_a$ or
$\ell_x$ distinct from the already existing ones cannot cross leaves
$\ell_a, \si_d^k(\ell_x), \dots, \si^{k(m-1)}(\ell_a),
\si_d^{km}(\ell_x)$ because either it would mean that leaves from
the same geolamination are linked, or it would contradict
Corollary~\ref{c:nothreelink}.}

\smallskip

By Claim A, we have $b_m\in (y_m, a_0]$. Consider possible locations of
$b_m$.

(1) If $x_0<b_m\le a_0$ then $\ol{a_mb_m}$ is linked with
$\ol{x_my_m},$ $\ol{x_{m+1}y_{m+1}}$ and $\ol{x_0y_0}$, which, by
Corollary~\ref{c:nothreelink}, implies that
$\ol{x_{m+1}y_{m+1}}=\ol{x_0y_0}$, and we are done (observe that, in
this case, by Lemma~\ref{l:sameperiod}, points $a_0, b_0$ are periodic
of the same period as $x_0, y_0$).

(2) The case $x_0=b_m$ is impossible because if $x_0=b_m$, then, by the order
preservation and by Claim A, the leaf
$\ol{x_{m+1}y_{m+1}}=\si_d^{k(m+1)}(\ell_x)$ is forced to be linked
with $\ell_a$, a contradiction.

(3) Otherwise we have $y_m<b_m<x_0$ and hence, by the order
preservation, $y_m\le x_{m+1}<b_m$. Then, by Claim A and because
images of $\ell_x$ do not cross, $b_m<y_{m+1}\le x_0$. Suppose that
$y_{m+1}=x_0$ while $y_0\ne x_1$. Applying $\si_d^k$ to leaves $\ol{x_{m+1}x_0}$ and
$\ol{x_0y_0}$ and using Claim A we see that $y_0\le x_{m+2}<x_1$.
However, the order preservation then implies that
$\ol{a_{m+1}b_{m+1}}$ crosses both $\ol{x_{m+1}x_0}$ and
$\ol{x_{m+2}x_1}$ and therefore crosses $\ell_a$ itself, a
contradiction. Hence the situation when $y_{m+1}$ coincides with $x_0$
can only happen if $y_0=x_1$. It follows that
then $\si^k_d(\ol{x_{m+1}y_{m+1}})=\ol{x_0y_0}$,
and we are done (as before, we need to rely on
Lemma~\ref{l:sameperiod} here).

Otherwise $b_m<y_{m+1}<x_0$ and the
arguments can be repeated as leaves $\si_d^{ki}(\ell_x), i=0,
\dots, m+1$ are circularly ordered. Thus, either $\ell_x$ is
periodic, $\ol{x_ny_n}=\ol{x_0y_0}$ for some $n$, and all leaves in
the $\si^k_d$-orbit of $\ell_x$ are circularly ordered, or the leaves
$\ol{x_iy_i}$ converge monotonically to a point of $\uc$. The latter
is impossible since $\si_d^k$ is expanding.  By
Lemma~\ref{l:sameperiod}, the leaf $\ell_a$ is periodic and its endpoints have
the same period as the endpoints of $\ell_x$.
\end{proof}

Theorem~\ref{t:compgap} is the main result of this section.

\begin{thm}\label{t:compgap}
Consider linked chords $\ell_a$, $\ell_x$
with mutually order preserving accordions, and set $B=\ch(\ell_a,\ell_x)$.
Suppose that not all forward images of $B$ have pairwise disjoint interiors.
Then there exists a finite periodic stand alone gap $Q$ such that
all vertices of $Q$ are in the forward orbit of $\si^r_d(B)$
for some minimal $r$, they belong to two, three, or four distinct periodic orbits
of the same period,
and the remap of $Q\cap\uc$ is not the identity unless
$Q=\si^r_d(B)$ is a quadrilateral.
\end{thm}

\begin{proof}
We may assume that there are two forward images of $B$ with
non-disjoint interiors. Choose the least $r$ such that the interior
of $\si^r_d(B)$ intersects some forward images of $B$. We may assume
that $r=0$ and, for some (minimal) $k>0$, the interior of the set
$\si^k_d(B)$ intersects the interior of $B$ so that
$\si_d^k(\ell_x)\in A_{\ell_x}(\ell_a)$. We write $x_i$, $y_i$ for
the endpoints of $\si_d^{ik}(\ell_x)$, and $a_i$, $b_i$ for the
endpoints of $\si_d^{ik}(\ell_a)$. By Lemma~\ref{l:linkstruct}
applied to both leaves, by the assumption of mutually order
preserving accordions, and because leaves in the forward orbits of
$\ell_a, \ell_x$ are pairwise unlinked, we may assume without loss
of generality that, for some $m\ge 1$,
$$x_0<a_0<y_0\le x_1<b_0\le a_1< \dots \le x_m<b_{m-1}\le a_m<y_m<b_m$$
and $x_m=x_0$, $y_m=y_0$, $a_m=a_0$, $b_m=b_0$, i.e., we have the
situation shown in Figure \ref{f:linkstruct}. Thus, for every $i=0$,
$\dots$, $k-1$, there is a loop $L_i$ of alternately linked
$\si_d^k$-images of $\si_d^i(\ell_a)$ and $\si_d^i(\ell_x)$. If the
$\si_d^k$-images of $\si_d^i(\ell_a)$ are concatenated to each
other, then their endpoints belong to the same periodic orbit,
otherwise they belong to two distinct periodic orbits. A similar
claim holds for $\si_d^k$-images of $\si_d^i(\ell_x)$. Thus, the
endpoints of $B$ belong to two, three or four distinct periodic
orbits of the same  period (the latter follows by
Corollary~\ref{c:nothreelink} and Lemma~\ref{l:linkstruct}).
Set $\ch(L_i)=T_i$ and consider some cases.

\begin{figure}
\includegraphics[width=4.5cm]{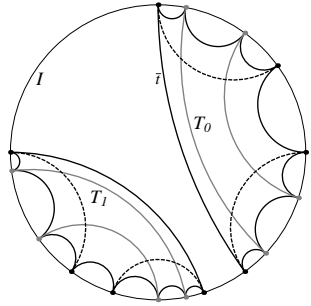}
{\caption{This figure illustrates
Theorem~\ref{t:compgap}(b) in the case $m>1$.}\label{f:compgap1}}
\end{figure}

(1) Let $m>1$ (this includes the ``flipping'' case from part (2) of
Corollary~\ref{c:nothreelink}). Let us show that the sets $T_i$
either coincide or are disjoint. Every image $\hell$ of $\ell_a$ in
$L_i$ crosses two images of $\ell_x$ in $L_i$ (if $m=2$ and $\ell_x$
is ``flipped'' by $\si^k_d$, we still consider $\ell_x$ and
$\si^k_d(\ell_x)$ as distinct leaves). By
Corollary~\ref{c:nothreelink}, no other image of $\ell_x$ crosses
$\hell$.

Suppose that interiors of $T_i$ and $T_j$ intersect. Let $\bt$ be an
edge of $T_i$ and $I=H_{T_i}(\bt)$ be the corresponding hole of
$T_i$. Then the union of two or three images of $\ell_a$ or $\ell_x$
from $L_i$ \emph{separates} $I$ from $\uc\sm I$ in $\cdisk$ (meaning
that any curve connecting $I$ with $\uc\sm I$ must intersect the
union of these two or three images of $\ell_a$ or $\ell_x$, see
Figure \ref{f:compgap1}). Hence if there are vertices of $T_j$
in $I$ \emph{and} in $\uc\sm I$ then there is a leaf of $L_j$
crossing leaves of $L_i$, a contradiction with the above and
Corollary~\ref{c:nothreelink}. Thus, the only way $T_i\ne T_j$ can
intersect is if they share a vertex or an edge. We claim that
this is impossible. Indeed, $T_i\ne T_j$ cannot share a vertex as
otherwise this vertex must be $\si_d^k$-invariant while all vertices of any
$T_r$ map to other vertices (sets $T_r$ ``rotate'' under $\si_d^k$).
Finally, if $T_i$ and $T_j$ share an edge $\ell$ then the same
argument shows that $\si_d^k$ cannot fix the endpoints of $\ell$,
hence it ``flips'' under $\si_d^k$. However this is impossible as
each set $T_r$ has at least four vertices and its edges ``rotate''
under $\si_d^k$.

So, the component $Q_i$ of $X=\bigcup^{k-1}_{i=0} T_i$ containing
$\si_d^i(\ell_a)$ is $T_i$. By Lemma~\ref{l:linkstruct}, the map
$\si_d|_{T_i\cap\uc}$ is order preserving or reversing. As
$\si_d$ preserves order on any single accordion,
$\si_d|_{T_i\cap\uc}$ is order preserving. The result now follows;
note that the first return map on $Q$ is not the identity map.

(2) Let $m=1$. This corresponds to part (3) of
Corollary~\ref{c:nothreelink}: both $\ell_a$ and $\ell_x$ have
endpoints of minimal period $k$, and the orbit of $\ell_a$ ($\ell_x$)
consists of $k$ pairwise disjoint leaves. Note that $T_0$ is a
quadrilateral, and the first return map on $T_0$ is the identity
map. Consider the case when not all sets $T_i$ are pairwise
disjoint. Note that, by the above, $T_0$ is a periodic stand alone
gap satisfying the assumptions of Proposition \ref{p:forconcat}. It
follows that every component of the union of $T_i$ is a
concatenation of gaps sharing edges with the same polygon. See
Figure \ref{f:compgap2}, in which the polygon is a triangle.
\begin{figure}
\includegraphics[width=4.5cm]{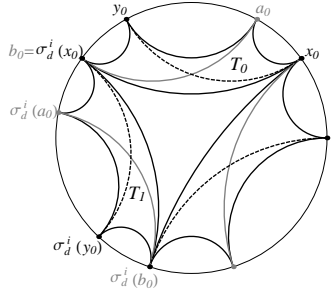}
{\caption{This figure illustrates Theorem~\ref{t:compgap}(b)
in the case $m=1$.}\label{f:compgap2}}
\end{figure}
\end{proof}

For a leaf $\ell_1\in \lam_1$, let $\B_{\lam_2}(\ell_1)$ be the
collection of all leaves $\ell_2\in \lam_2$ which are linked with
$\ell_1$ and have mutually order preserving accordions with
$\ell_1$. Observe that if $\ell_1$ is (pre)critical, then
$\B_{\lam_2}(\ell_1)=\0$ by Definition~\ref{d:opacc}. Similarly, no
leaf from $\B_{\lam_2}(\ell_1)$ is (pre)critical.

\begin{cor}\label{c:finaccord}
The collection $\B_{\lam_2}(\ell_1)$ is finite.
\end{cor}

\begin{proof}
Suppose first that $\ell_1$ is not (pre)periodic. Let us show that
the convex hull $B$ of $\ell_1$ and leaves $\n_1, \dots, \n_s$ from
$\B_{\lam_2}(\ell_1)$ is wandering. By Theorem~\ref{t:compgap}, for
each $i$, the set $B_i=\ch(\ell_1, \n_i)$ is wandering (because
$\ell_1$ is not (pre)periodic). This implies that if $i\ne j$ then
$\si^i_d(\ell_1)$ and $\si_d^j(\n_t)$ are disjoint (otherwise
$\si_d^i(B_t)$ and $\si_d^j(B_t)$ are non-disjoint). Moreover,
$\si_d^i(\ell_1)$ and $\si_d^j(\ell_1)$ are disjoint as otherwise, by
Lemma~\ref{l:conconv}, the leaf $\ell_1$ is (pre)periodic. Therefore
$\si_d^j(\ell_1)$ is disjoint from $\si_d^i(B)$.

Suppose that $\si_d^i(B)$ and $\si_d^j(B)$ are non-disjoint. By
the just proven then, say, $\si_d^j(\n_1)$ is
non-disjoint from $\si_d^i(B)$. Again by the just proven
$\si_d^j(\n_1)$ is disjoint from $\si_d^i(\ell_1)$. Hence the only
possible intersection is between $\si_d^j(\n_1)$ and, say,
$\si_d^i(\n_2)$. Moreover, since $\si_d^j(\ell_1)$ is disjoint from
$\si_d^i(B)$, then $\si_d^j(\n_1)\ne \si_d^i(\n_2)$ and, moreover,
as distinct leaves of the same lamination, the leaves
$\si_d^j(\n_1), \si_d^i(\n_2)$ cannot cross. Hence the only way
$\si_d^j(\n_1)$ and $\si_d^i(\n_2)$ are non-disjoint is that
$\si_d^j(\n_1)$ and $\si_d^i(\n_2)$ are concatenated.

Assume that $\si^t_d(\n_2)$ is concatenated with
$\n_1$ at an endpoint $x$ of $\n_1$. Clearly, $x$ is a common vertex
of $B$ and of $\si^t_d(B)$. Hence $\si^t_d(x)$ is a common vertex of
$\si^t_d(B)$ and $\si^{2t}_d(B)$, etc. Connect points $x$,
$\si^t(x),$ $\si^{2t}(x),$ $\dots$ with consecutive chords $\m_0,$
$\m_1,$ $\dots$. These chords are pairwise unlinked because, as it
follows from the above, the sets $\si_d^r(B)$, $r=0$, $1$, $\dots$ have pairwise
disjoint interiors. Hence, by Lemma~\ref{l:concat}, the point $x$ is
(pre)periodic, a contradiction with the fact that all sets
$B_i=\ch(\ell_1, \n_i)$ are wandering. Thus, $B$ is wandering. Hence,
by \cite{kiw02}, the collection $\B_{\lam_2}(\ell_1)$ is finite.

Suppose now that $\ell_1$ is periodic. Then by
Theorem~\ref{t:compgap} any leaf of $\B_{\lam_2}(\ell_1)$ is
periodic with the same periods of endpoints. This implies that in
this case the collection $\B_{\lam_2}(\ell_1)$ is finite. Finally,
if $k>0$ is the minimal number such that $\si_d^k(\ell_1)$ is
periodic and $\ell_2\in \B_{\lam_2}(\ell_1)$ then $\si_d^k(\ell_2)$
is linked with $\si_d^k(\ell_1)$ which implies that $\ell_2$ is a
$\si_d^k$-preimage of one of finitely many leaves from
$\B_{\lam_2}(\si_d^k(\ell_1))$. Thus, in this case
$\B_{\lam_2}(\ell_1)$ is finite too.
\end{proof}

\section{Linked Quadratically Critical Geolaminations}\label{s:qcrit}

The main result of Section \ref{s:qcrit} is that two linked or essentially
coinciding geolaminations with qc-portraits have the same perfect
sublamination (see Definition~\ref{d:qclink1}). In this section, we will always assume that the
laminations $(\lam_1,\qcp_1)$ and $(\lam_2,\qcp_2)$ are linked or essentially equal.

\subsection{Smart Criticality}\label{ss:smart}

Our aim in Subsection~\ref{ss:smart} is to introduce \emph{smart criticality},
a principle which allows one to use a flexible choice of critical
chords of $\lam_1$ and $\lam_2$ in order to treat certain sets of
linked leaves of $\lam_1$ and $\lam_2$ as if they were sets of one
lamination.

\begin{lem}\label{l:accorder}
If $\ell_1\in\lam_1$ is not a special critical leaf, then each critical set $C$
of $\qcp_2$ has a spike $\ol{c}$ unlinked with $\ell_1$; these spikes
form a full collection $\mathcal E$ of spikes of $\lam_2$ unlinked with
$\ell_1$. If
an endpoint $x$ of $\ell_1$ is neither a vertex of a special critical
cluster nor a common vertex of associated critical quadrilaterals of our
geolaminations, then $\mathcal E$ can
be chosen so that $x$ is not an endpoint of a spike from $\mathcal E$.
\end{lem}

\begin{proof}
Since $\ell_1$ is not a special critical leaf, spikes of $\lam_2$ from special
critical clusters are unlinked with $\ell_1$. Otherwise take a pair of
associated critical \ql s $A\in \lam_1, B\in \lam_2$ with non-strictly
alternating on $\uc$ vertices

$$a_0\le b_0\le a_1\le b_1\le a_2\le b_2\le a_3\le b_3\le a_0$$

\noindent and observe, that $\ell_1$ is contained, say, in $[a_0,
a_1]$ and hence is unlinked with the spike $\ol{b_1b_3}$ of $B$. The
last claim is left to the reader.
\end{proof}

Denote by $\mathcal E_{\lam_2}(\ell_1)$ a full collection of spikes from Lemma~\ref{l:accorder}.


\begin{cor}\label{c:accorder}
If $\ell_1=\ol{ab}\in\lam_1$ is not a special critical leaf, then
$A=A_{\lam_2}(\ell_1)$ is contained in the closure of a component of
$\disk\sm \mathcal E_{\lam_2}(\ell_1)^+$, 
and $\si_d|_{A\cap\uc}$ is $($non-strictly$)$ monotone. Let
$\ell_2=\ol{xy}\in \lam_2$ and $\ell_1\cap \ell_2\ne \0$. Then:

\begin{enumerate}

\item if $\ell_1$ and $\ell_2$ are concatenated at a point $x$
    that is neither a vertex of a special critical cluster nor
    a common vertex of associated critical quadrilaterals of our
    geolaminations, then $\si_d$ is (non-strictly) monotone on $\ell_1\cup \ell_2$;

\item if $\ell_2$ crosses $\ell_1$, then, for each $i$, we have
$\si_d^i(\ell_1)\cap\si_d^i(\ell_2)\ne \0$, and one the following
holds:

\begin{enumerate}

\item $\si_d^i(\ell_1)=\si_d^i(\ell_2)$ is a point or a leaf shared by $\lam_1, \lam_2$;

\item $\si_d^i(\ell_1), \si_d^i(\ell_2)$ share an endpoint;

\item $\si_d^i(\ell_1), \si_d^i(\ell_2)$ are linked and have
the same order of endpoints as $\ell_1, \ell_2$;

\end{enumerate}

\item points $a$, $b$, $x$, $y$
are either all $($pre$)$\-pe\-ri\-odic of the same eventual period, or are
all not $($pre$)$periodic.

\end{enumerate}

\end{cor}

\begin{proof}
Set $\mathcal E=\mathcal E_{\lam_2}(\ell_1)$. 
If $\ell_1$ coincides with one of spikes from $\mathcal E$, the claim follows
(observe that then by definition $A=\ell_1$ as spikes of sets of
$\lam_2$ do not cross leaves of $\lam_2$). Otherwise there exists a
unique complementary component $Y$ of $\mathcal E^+$ with
$\ell_1\subset Y$ (except perhaps for the endpoints). The fact that
each leaf of $\lam_2$ is unlinked with spikes from $\mathcal E$
implies that $A_{\lam_2}(\ell_1)\subset \ol{Y}$. This proves the
main claim of the lemma.

(1) By Lemma~\ref{l:accorder}, the collection $\mathcal E$ can be
chosen so that $x$ is not an endpoint of a chord from $\mathcal E$.
The construction of $Y$ then implies that $\si_d$ is monotone on
$\ell_1\cup \ell_2$.

(2) We use induction. By Definition~\ref{d:qclink1}, if a critical leaf
$\n_1\in \lam_1$ crosses a leaf $\m_2\in \lam_2$ and comes from a
special critical cluster, then both $\n_1$ and $\m_2$ come from a
special critical cluster and have the same image. Thus we may assume
that neither $\si^i_d(\ell_1)$ nor $\si^i_d(\ell_2)$ are from a special
critical cluster. We may also assume that $\si^i_d(\ell_1)$ and
$\si^i_d(\ell_2)$ do not share an endpoint as otherwise the claim is
obvious. Hence it remains to consider the case when $\si^i_d(\ell_1)$
and $\si^i_d(\ell_2)$ are linked and are not special critical leaves. Then by
the main claim either their images are linked or at least they share an
endpoint.

(3) By Lemma~\ref{l:conconv}, if an endpoint of a leaf of a
geolamination is (pre)\-pe\-ri\-odic, then so is the other endpoint
of the leaf. Consider two cases. Suppose first that an image of
$\ell_1$ and an image of $\ell_2$ ``collide'' (i.e., have a common
endpoint $z$). By the above, if $z$ is (pre)periodic, then all
endpoints of our leaves are, and if $z$ is not (pre)periodic, then
all endpoints of our leaves are not (pre)periodic. Suppose now that
no two images of $\ell_1, \ell_2$ collide. Then it follows that
$\ell_1$ and $\ell_2$ have mutually order preserving accordions, and
the claim follows from Theorem~\ref{t:compgap}.
\end{proof}

Lemma~\ref{l:accorder} and Corollary~\ref{c:accorder} implement
smart criticality. Indeed, given a geolamination $\lam$, a gap or leaf $G$
of it is such that the set $G\cap \uc$ (loosely) consists of points whose orbits avoid critical
sets of $\lam$. It follows that any power of the map is order
preserving on $G\cap \uc$. It turns out that we can treat sets $X$
formed by linked leaves of two linked/essentially equal
geolaminations similarly by varying our choice of the full
collection of spikes on each step so that the orbit of $X$ avoids
\emph{that particular} full collection of spikes on \emph{that
particular} step (thus \emph{smart} criticality). Therefore,
similarly to the case of one geolamination, any power of the map is
order preserving on $X$. This allows one to treat such sets $X$
almost as sets of one geolamination.

Lemma~\ref{l:accorder2} describes how $\si_d$ can be
\emph{non-strictly} monotone on $A\cap\uc$ taken from
Corollary~\ref{c:accorder}. A concatenation $\rc$ of spikes of a
geolamination $\lam$ such that the endpoints of its chords are
monotonically ordered on the circle will be called a \emph{chain of
spikes $($of $\lam)$}.

\begin{lem}\label{l:accorder2}
Suppose that $\ell_a=\ol{ab}\in \lam_1$, $\ell_x=\ol{xy}\in \lam_2$,
$a<x<b\le y<a$ $($see Figure \ref{f:accorder2}$)$ and if $b=y$, then
$b$ is neither a vertex of a special critical cluster nor a common
vertex of associated critical quadrilaterals of our geolaminations. Let
$\si_d(a)=\si_d(x)$. Then either both $\ell_a$, $\ell_x$ are contained
inside the same special critical cluster, or there are chains of spikes $\rc_1$ of $\lam_1$
and $\rc_2$ of $\lam_2$ connecting $a$ with $x$. If one of leaves
$\ell_a$, $\ell_x$ is not critical, we may assume that
$\rc_1^+\cap\uc\subset [a, x]$ and that $\rc_2^+\cap\uc\subset [a, x]$.
\end{lem}

Recall that, according to our terminology, a chord is contained \emph{inside}
a special critical cluster $S$ if it is a subset of $S$ intersecting
the interior of $S$.

\begin{proof}
First assume that one of the leaves $\ell_a, \ell_x$ (say, $\ell_a$) is
a special critical leaf. Then both $a$ and $b$ are vertices of a special
critical cluster. By the assumptions, this implies that $b\ne y$ and
hence $\ell_a$ and $\ell_x$ are linked and are inside
a special critical cluster. Assume from now on that
neither $\ell_a$ nor $\ell_x$ is a special critical leaf.

\begin{figure}
\includegraphics[width=4.5cm]{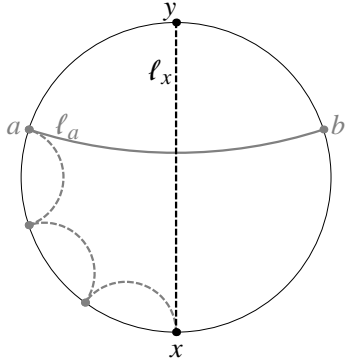}
{\caption{This figure illustrates
Lemma~\ref{l:accorder2}. Here the leaves $\ell_a,\ell_x$
collapse around a chain of spikes shown as dashed grey geodesics.}\label{f:accorder2}}
\end{figure}

By Lemma~\ref{l:accorder}, choose a full collection $\mathcal A_2$
of spikes of $\lam_2$ unlinked with $\ell_a$ and a full collection
$\mathcal A_1$ of spikes of $\lam_1$ unlinked with $\ell_x$. By the
assumptions and Lemma~\ref{l:accorder}, we may choose these
collections so that if $b=y$, then $b=y\nin \mathcal A_1^+\cup
\mathcal A_2^+$. Thus in any case the point $\ell_a\cap \ell_x=w\in \disk$
does not belong to $\mathcal A_1^+\cup \mathcal A_2^+$.

It follows that there is a well-defined component $Y$ of $\cdisk\sm
\mathcal A_i^+$ containing $\ell_a\cup \ell_x$ except perhaps for
the endpoints. Since $\si_d(a)=\si_d(x)$, there is a chain of spikes
$\rc_2\subset \mathcal A_2$ of $\lam_2$ and a chain of spikes $\rc_1\subset \mathcal A_1$
of $\lam_1$ connecting $a$ and $x$. Suppose that, say,
$\rc_1^+\cap\uc\subset [x, a]$. Since all spikes are critical chords
which cross neither $\ell_a$ nor $\ell_x$, this implies that both
$\ell_a$ and $\ell_x$ are critical. Therefore, if at least one of the
leaves $\ell_a$, $\ell_x$ is not critical, then we may assume that
$\rc_1^+\cap\uc\subset [a, x]$ and that $\rc_2^+\cap\uc\subset [a, x]$.
\end{proof}

The assumptions of Lemma~\ref{l:accorder2} automatically hold if leaves
$\ell_a$, $\ell_x$ are linked and one of them (say, $\ell_a$) is
critical; in this case, by Corollary~\ref{c:accorder},
the point $\si_d(\ell_a)$ is
an endpoint of $\si_d(\ell_x)$, and, renaming the points, we may assume
that $\si_d(a)=\si_d(x)$.

\begin{dfn}\label{d:accorder2}
Non-disjoint leaves $\ell_1\ne \ell_2$ are said to \emph{collapse
around chains of spikes} if there are two chains of spikes, one in each of the two
geolaminations, connecting two adjacent endpoints of $\ell_1$, $\ell_2$ as
in Lemma~\ref{l:accorder2}.
\end{dfn}

Smart criticality allows one to treat accordions as gaps of one
geolamination provided images of leaves do not collapse around chains of spikes.

\begin{lem}\label{l:indepcrit}
Let $\ell_1$, $\ell_2$ be linked leaves from $\lam_1$, $\lam_2$ such that
there is no $t$ with $\si_d^t(\ell_1)$, $\si_d^t(\ell_2)$ collapsing
around chains of spikes. Then there exists an $N$ such that the
$\si_d^N$-images of $\ell_1$, $\ell_2$ are linked and have mutually
order preserving accordions. Conclusions of Theorem~\ref{t:compgap}
hold for $\ell_1$, $\ell_2$, and $B=\ch(\ell_1, \ell_2)$ is either
wandering or $($pre$)$periodic so that $\ell_1$, $\ell_2$ are
$($pre$)$periodic of the same eventual period of endpoints.
\end{lem}

\begin{proof}
By way of contradiction, suppose that there exists the minimal $t$
such that $\si^{t+1}_d(\ell_1)$ is not linked with
$\si^{t+1}_d(\ell_2)$. Then $\si_d^t(\ell_1)$ crosses
$\si_d^t(\ell_2)$ while their images have a common endpoint. Hence
Lemma~\ref{l:accorder2}, applied to $\si_d^t(\ell_1)$ and
$\si_d^t(\ell_2)$, implies that $\si_d^t(\ell_1)$, $\si_d^t(\ell_2)$
collapse around a chain of spikes, a contradiction. Thus,
$\si_d^t(\ell_1)$ and $\si_d^t(\ell_2)$ cross for any $t\ge 0$. In
particular, no image of either $\ell_1$ or $\ell_2$ is ever
critical.

By Lemma~\ref{l:conconv}, choose $N$ so that leaves
$\si_d^N(\ell_1)=\ol{ab}$ and $\si_d^N(\ell_2)=\ol{xy}$ are periodic
or have no (pre)\-pe\-ri\-odic endpoints. If $\ol{ab}$ and $\ol{xy}$
are periodic, then no collapse around chains of critical leaves on
\emph{any} images of $\ol{ab}, \ol{xy}$ is possible (for
set-theoretic reasons). Hence $\si_d^N(\ell_1), \si_d^N(\ell_2)$ are
linked and have mutually order preserving accordions as desired.

Suppose now that our leaves have non-(pre)\-pe\-ri\-odic endpoints.
Evidently, the set $E$ of all endpoints of all possible chains of
spikes is finite. Thus, there exists an $N$ such that if $n\ge N$,
then $\si_d^n(a)$ is disjoint from $E$.
The same
holds for $b$, $x$ and $y$, so we may assume that, for $n\ge N$, no
endpoint of $\si_d^n(\ell_1)$ or $\si_d^n(\ell_2)$ is in  $E$. Hence,
the $\si_d^N$-images of $\ell_1$, $\ell_2$ are linked
and have mutually order preserving accordions.
\end{proof}

\subsection{Linked perfect laminations}\label{ss:consmart}

\begin{lem}\label{l:countable}
The set $\Tc$ of all leaves of $\lam_2$ non-disjoint from a leaf
$\ell_1$ of $\lam_1$ is at most countable. Thus, if $\ell_1$ is an
accumulation point of uncountably many leaves of $\lam_1$ then
$\ell$ is unlinked with any leaf of $\lam_2.$
\end{lem}

\begin{proof}
If $\ell_1$ has (pre)periodic endpoints, then, by
Corollary~\ref{c:accorder}, any leaf of $\lam_2$ non-disjoint from
$\ell_1$ has (pre)periodic endpoints implying the first claim
of the lemma in this case. Let $\ell_1$ have no (pre)periodic
endpoints. Then, by Corollary~\ref{c:accorder}, leaves of $\lam_2$
non-disjoint from $\ell_1$ have no (pre)\-pe\-ri\-odic endpoints. By
Lemma~\ref{l:cones}, there are finitely many leaves with an endpoint
being a given eventual image of an endpoint of $\ell_1$. Hence
the set of all leaves of $\lam_2$ with an endpoint
whose orbit collides with the orbit of an endpoint of $\ell_1$ is
countable. If we remove them from $\Tc$, we will get a new
collection $\Tc'$ of leaves $\ell_2'$, which have mutually order
preserving accordions with $\ell_1$. By Corollary~\ref{c:finaccord}, the collection
$\Tc'$ is finite. This completes the proof of the first claim of the
lemma.
The second claim follows immediately. \end{proof}

Let $\qcp$ be a qc-portrait of a geolamination $\lam$. Since, by
Lemma~\ref{l:cridisj}, distinct critical sets of the perfect
sublamination $\lam^c$ are disjoint, each critical set of $\lam$ is
contained in a unique critical set of $\lam^c$. Hence $\qcp$ generates
the \emph{critical pattern $\zc(\qcp)$ of $\qcp$ in $\lam^c$}, and so
each geolamination with critical portrait $(\lam, \qcp)$ gives rise to
the perfect geolamination with critical pattern $(\lam^c, \zc(\qcp))$.

\begin{thm}\label{t:noesli}
If $(\lam_1, \qcp_1)$ and $(\lam_2, \qcp_2)$ are geolaminations
with qc-portraits that are linked or essentially equal, then we
have the following equality:
$(\lam^c_1, \zc(\qcp_1))$ $=$ $(\lam^c_2, \zc(\qcp_2))$.
\end{thm}

\begin{proof} By way of contradiction, assume that
$\lam^c_1\not\subset \lam^c_2$; then
$\lam^c_1\not\subset \lam_2$, and there exists a leaf $\ell^c_1\in
\lam^c_1\sm \lam_2$. Then, by Lemma~\ref{l:countable}, the leaf
$\ell^c_1$ is inside a gap $G$ of $\lam_2$.
Since $\lam^c_1$ is perfect, from at least
one side all one-sided neighborhoods of $\ell^c_1$ contain
uncountably many leaves of $\lam^c_1$. Hence $G$ is uncountable (if
$G$ is finite or countable, then there must exist edges of $G$ which
cross leaves of $\lam^c_1$, a contradiction as above). Thus, there
are uncountably many leaves of $\lam^c_1$ inside $G$; these leaves
connect points of $G\cap\uc$. This contradicts Corollary~\ref{c:noinf}.
\end{proof}

Jan Kiwi showed in \cite{kiwi97} that if all critical sets of a geolamination
$\lam$ are critical
leaves with \emph{aperiodic kneading}, then its perfect sublamination
$\lam^c$ is completely determined by these critical leaves (he also showed that
this defines the corresponding lamination $\sim$ such that $\lam^c=\lam_\sim$
and that $\sim$ is dendritic).
Our results are related to Kiwi's. Indeed, by Theorem~\ref{t:noesli},
if $\lam$ is a geolamination with a qc-portrait $\qcp$, then
$\lam^c\subset \lam$ is completely
defined by $\qcp$; in other words, if there is another geolamination
$\hlam$ with the same qc-portrait $\qcp$, then still $\hlam^c=\lam^c$.
However, Theorem~\ref{t:noesli} takes
the issue of how critical data impacts the perfect sublamination of a
geolamination further as it considers the dependence of the perfect
sublaminations upon critical data while relaxing the conditions on
critical sets and allowing for ``linked perturbation'' of the critical
data. Therefore, Theorem~\ref{t:noesli} could be viewed as a rigidity
result: ``linked perturbation'' of critical data does not change the
perfect geolamination.

\begin{dfn}\label{d:linoreq}
Let $\lam_1$ and $\lam_2$ be geolaminations. Suppose that there are
geolaminations with qc-portraits $(\lamm_1, \qcp_1)$, $(\lamm_2,
\qcp_2)$ such that $\lam_1\subset \lamm_1,$ $\lam_2\subset \lamm_2$
and $(\lamm_1, \qcp_1)$ and $(\lamm_2, \qcp_2)$ are linked
(essentially equal). Then we say that $\lam_1$ and $\lam_2$ are
\emph{linked} (\emph{essentially equal}, respectively).
\end{dfn}

Theorem~\ref{t:noesli} immediately implies Lemma~\ref{l:disjoint}.

\begin{lem}\label{l:disjoint}
Let $\lam_1$ and $\lam_2$ be geolaminations that are linked or
 essentially equal and such that the
geolaminations $\lamm_1$ and $\lamm_2$ from Definition~\ref{d:linoreq}
have perfect sublaminations equal to perfect sublaminations
$\lam^c_1\subset \lam_1$ and $\lam^c_2\subset \lam_2$. Then
$\lam^c_1=\lam^c_2=\lam^c$, and critical patterns of\, $\qcp_1$ in
$\lam^c$ and of\, $\qcp_2$ in $\lam^c$ coincide.
\end{lem}

The second condition above means that by inserting (if necessary)
critical \ql s into critical sets of $\lam_1$ and $\lam_2$ we do not
change the perfect sublamination of either geolamination.

\begin{cor}\label{c:disjoint}
Let $\lam_1$ and $\lam_2$ be geolaminations that are linked or
essentially equal.
 Suppose that all critical sets
in both $\lam_1$ and $\lam_2$ are finite. Then
$\lam^c_1=\lam^c_2=\lam^c$, and critical patterns of $\qcp_1$ in
$\lam^c$ and of $\qcp_2$ in $\lam^c$ coincide.
\end{cor}

\begin{proof}
Choose geolaminations $\lamm_1$ and $\lamm_2$ from
Definition~\ref{d:linoreq}. These are constructed by inserting (if
necessary) \ql s into critical sets of $\lam_1$ and $\lam_2$ and then mapping
them forward and pulling them back. Since the critical sets of $\lam_1$, $\lam_2$
are finite, this creates no new non-isolated leaves (such leaves can
only be created if the grand orbits of inserted \ql s accumulate inside
infinite critical gaps). Hence the perfect sublamination of $\lamm_1$
equals $\lam^c_1$ and the perfect sublamination of $\lamm_2$ equals
$\lam^c_2$. By Lemma~\ref{l:disjoint} this implies the desired.
\end{proof}

\section{Applications}\label{l:appli}
A polynomial is \emph{dendritic} if its Julia set is
connected, and all its periodic
points are repelling. The main theorem of this section
gives a combinatorial model for the space $\md_3$ of all
cubic critically marked dendritic polynomials $P$.
Recall that, by
Kiwi \cite{kiwi97}, if $P$ is a dendritic polynomial, then $P|_{J(P)}$
is monotonically semiconjugate by a map $\Psi_P$ to its topological
polynomial $f_{\sim_P}$ on its topological Julia set $J_{\sim_P}$,
where $J_{\sim_P}$ is a dendrite all of whose points have
finite order. The lamination $\sim_P$ and the associated
geolamination $\lam_{\sim_P}$ are then called \emph{dendritic} (see
Definition~\ref{d:dendrilam}).

\subsection{Cubic dendritic geolaminations}\label{ss:cudebi}
Observe that all gaps of $\sim_P$ are finite. Recall that the \emph{degree} of
a gap $G$ of $\sim_P$ was defined right above Definition~\ref{d:cristuff}.
By a \emph{critical set} we mean either a
gap of degree greater than one or a critical leaf.

We consider the family $\ld_3$ of all cubic dendritic geolaminations $\lam$
with an ordered pair of critical sets and parameterize it with an
ordered pair of sets specifically related to those critical sets.
This gives a geometric interpretation of this family analogous to
Thurston's description of the geolamination $\qml$. For
a cubic geolamination $\lam$ with critical sets $C_1\ne C_2$
of degree two, we consider the \emph{co-critical set} of $C_1$ (i.e., the set with the
same image as $C_1$ but disjoint from $C_1$) and the \emph{minor}
set $\si_3(C_2)$. If $\lam$ has a unique critical set $C$ of degree
three, we associate $\lam$ with the pair of sets $(C,\si_3(C))$.

In Definition~\ref{d:unibi}, we mimic Milnor's terminology for
polynomials.

\begin{dfn}[Unicritical and bicritical geolaminations]\label{d:unibi}
A geolamination that has a critical set of degree three is called
\emph{unicritical}. Otherwise $\lam$ is said to be
\emph{bicritical}.
\end{dfn}

\emph{Full portraits} extend the notion of a critical pattern in the
cubic case. Observe that, in the definition below, we allow for the
possibility that two critical sets are non-disjoint (one could even
be an edge of the other one).

\begin{dfn}[Full portraits]\label{d:cripat3}
Consider a cubic geolamination $\lam$. An ordered pair $(C_1, C_2)$
of critical sets of $\lam$ is called a \emph{full portrait of
$\lam$} if either (1) $C=C_1=C_2$ is a unique critical set of a
unicritical geolamination $\lam$, or (2) $C_1\ne C_2$ and $\lam$
could be either bicritical or unicritical. Then the triple $(\lam,
C_1, C_2)$ is called a \emph{cubic geolamination with full
portrait}.
\end{dfn}

By Definition~\ref{d:cripat3}, if $\lam$ has two disjoint critical
sets $K_1$, $K_2$, then $(K_1, K_1)$ is \emph{not} a full portrait of
$\lam$. Thus, if $\lam$ is dendritic and bicritical, then a full
portrait of $\lam$ is just an ordering $(G,H)$ of the two critical
sets $G$ and $H$ of $\lam$.
In fact, all critical patterns of dendritic geolaminations are full
portraits. However, if $\lam$ has an all-critical triangle $T$, then
a pair ($T$, an edge of $T$) or a pair formed by two distinct edges
of $T$ are full portraits but are not critical patterns of $\lam$.
In general (not assuming that $\lam$ is dendritic), here are
possible cases for full portraits of $\lam$.
\begin{enumerate}
\item The geolamination $\lam$ has a unique critical set $C$
which is not an all-critical triangle. Then the only full portrait of $\lam$
is $(C,C)$.

\item The geolamination $\lam$ has a unique critical set, which is an all-critical
triangle $T$. 
Then full portraits of $\lam$ are $(T,T)$,
an ordered pair of edges of $T$, and $T$ with one of its edges.

\item The geolamination $\lam$ has two disjoint critical sets.
Then the two orderings of these sets are the only two full
portraits of $\lam$.

\item The geolamination $\lam$ has two non-disjoint critical sets of degree two. Some
cases of this type were considered in (2) above. Otherwise $\lam$
can have two critical sets sharing a vertex/a non-critical
edge. The critical sets of $\lam$ share a critical leaf
$\ovc$ only if one critical set is a critical leaf $\ovc$ and the
other critical set is a gap $G$ which has $\ovc$ as its edge and is
otherwise mapped two-to-one to its image.
\end{enumerate}

In most cases, if we choose two critical chords in two different
critical sets of a full portrait $\fup$, one in each critical set, then these two
critical chords are distinct and 
form a qc-portrait, whose
elements are contained in elements of $\fup$. However, there are
exceptions. Indeed, if $T$ is an all-critical triangle, then there
may exist full portraits, whose two sets share a critical chord.
Otherwise, in case (4), we may have two critical sets, one of which is
a critical gap $G$ with a critical edge $\ovc$, while the other one
is $\ovc$ itself.

In Definition~\ref{d:admicri}, we introduce the
\emph{co-critical} set $\coc(C)$ of a set $C$. Note that co-critical
sets are only defined for some gaps or leaves $C$ of a geolamination
and are in general not gaps or leaves of the same geolamination. If
a cubic geolamination $\lam$ has a unique critical set $C$ of degree three then
no hole of $C$ is greater than $\frac13$ while any other set has a
unique hole of length greater than $\frac13$. Suppose that $\lam$
has two critical sets $C_1, C_2$. Then either set is of degree two
and has a unique hole of length greater than $\frac13$. 
It
is easy to see that if a leaf or gap $D$ separates $C_1$ and $C_2$
then it has two holes of length greater than $\frac13$, otherwise
$D$ has a unique such hole.

\begin{dfn}[co-critical set]\label{d:admicri}
Let $C$ be a (possibly degenerate) leaf or a gap of a cubic geolamination 
$\lam$. Moreover, assume that either $C$ is critical of degree three, or
there is exactly one hole of $C$ of length at least $\frac13$.
If $C$ is of degree three, we set $\coc(C)=C$.
Otherwise let $H$ be a unique hole of $C$ of length $\ge\frac 13$.
Let $A$ denote the set of all points in $\ol H$ with the images in $\si_3(C)$.
Set $\coc(C)=\ch(A)$.
The set $\coc(C)$ is called the \emph{co-critical set} of $C$.
\end{dfn}

Definition~\ref{d:minor} mimics Thurston \cite{thu85}.

\begin{dfn}[minor set]\label{d:minor}
Let $(\lam, C, D)$ be a geolamination with full portrait. Then
$\si_3(D)$ is called the \emph{minor set of $(\lam, C, D)$}.
\end{dfn}

We are ready to define tags of cubic geolaminations with full portraits.

\begin{dfn}[mixed tag]\label{d:siblita} Suppose that $(\lam, \fup)$
is a cubic geolamination with full portrait $\fup=(C_1, C_2)$. Then we
call the set $\ta(\lam, \fup)=\ta(\fup)=\coc(C_1)\times
\si_3(C_2)\subset \cdisk\times \cdisk$
the \emph{mixed tag} of $(\lam, \fup)$.
\end{dfn}

It is easy to see that sets $\coc(C_1)$ (and hence
mixed tags) are well-defined. Note also that the mixed tag $T$ of a
cubic geolamination is the product of two sets, each of which is a
point, a leaf, or a gap. We can think of $T\subset \cdisk\times
\cdisk$ as a higher dimensional analog of a gap or a leaf of a
geolamination in $\ol{\disk}$. We show that the union of
tags of all dendritic cubic geolaminations with full portraits is a
(non-closed) ``geolamination'' in $\ol{\disk}\times\ol{\disk}$.

Recall that a critical \ql{} is called \emph{collapsing} if its
image is a non-degenerate leaf. The proof of
Lemma~\ref{l:linkedgivesco} is left to the reader.

\begin{lem} \label{l:linkedgivesco}
Suppose that $C$ is a collapsing quadrilateral or a critical leaf.
Then $C$ is the convex hull of the set of points
$[\coc(C)+\frac13]\cup [\coc(C)+\frac23]$ in $\uc$. Moreover, if
$C_i$ $(i=1,2)$ are collapsing strongly linked quadrilaterals then
$\coc(C_1)$ and $\coc(C_2)$ are linked leaves.
\end{lem}

Proposition~\ref{p:colocation} helps dealing with co-critical sets.

\begin{prop}\label{p:colocation}
Suppose that $(\lam,\fup)$ is a cubic geolamination with a full
portrait, $C\in\fup$, and $\ell=\ol{ab}$ is an edge of $\coc(C)$
with $(a,b)\cap C=\0$. Then:
\begin{enumerate}
\item $\si_3|_{(a,b)}$ is one-to-one, and
\item if $D$ is the other element of $\fup$, then
$\si_3(D)\subset [\si_3(b),\si_3(a)]$.
\end{enumerate}
\end{prop}

\begin{proof}
If $(a, b)$ is of length $\frac13$ then there is nothing to prove.
If $(a,b)$ had length greater than $\frac13$, then there would be a
sibling of $\ell$ disjoint from $\ell$ with endpoints in $(a,b)$.
Evidently, such a leaf would be an edge of $C$, contradicting the
choice of $(a,b)$. Thus we may assume that the length of $(a, b)$ is
less than $\frac13$. This implies (1).

If $C=D$ is of degree three, then $\ell$ is an edge of $C$ and (2)
follows immediately. Otherwise, let $a'=a+\frac13$ and
$b'=b+\frac23$. Then $\ol{a'b'}\subset C$ and, for geometric reasons,
$D$ must have endpoints in $[b',a]\cup[b,a']$. Since each of these
intervals maps onto $[\si_3(b),\si_3(a)]$ one-to-one, (2) follows.
\end{proof}

Recall that strongly linked \ql s are defined in
Definition~\ref{d:strolin}. In particular, two strongly linked \ql s
may have common vertices.

\begin{prop}\label{p:linkco}
Suppose that $\ell_1$, $\ell_2$ are linked chords of\, $\uc$, whose
endpoints are contained in an interval of length at most $\frac13$.
Then $\coc(\ell_1)$ and $\coc(\ell_2)$ are strongly
linked collapsing quadrilaterals.
\end{prop}

\begin{figure}
 \includegraphics[width=4.5cm]{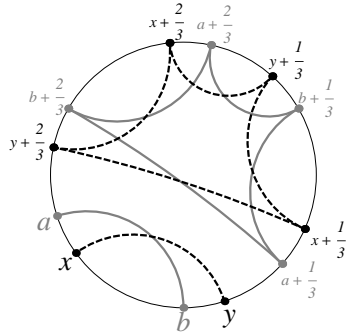}
 \caption{This figure illustrates Proposition~\ref{p:linkco}.  The linked
 chords $\ol{ab}$, $\ol{xy}$ lying in an interval of length less
 than $\frac13$ have $\coc({\ol{ab}})$, $\coc({\ol{xy}})$ as strongly
 linked quadrilaterals.} \label{f:cocri}
\end{figure}

\begin{proof}
As neither leaf is critical, we may assume $\ell_1=\ol{ab},
\ell_2=\ol{xy}$ and $a<x<b<y\le a+\frac13$.
Proposition~\ref{p:colocation} implies that vertices of \ql s
$\coc(\ell_1)$ and $\coc(\ell_2)$ satisfy the following
inequalities: $a+\frac13<x+\frac13<b+\frac13<y+\frac13\le
a+\frac23<x+\frac23<b+\frac23<y+\frac23\le a$, see Figure
\ref{f:cocri}.
\end{proof}

Recall that, by Definition~\ref{d:linoreq}, two geolaminations are
said to be \emph{linked} (\emph{essentially equal}) if they can be
``tuned'' to geolaminations with qc-portraits which indeed \emph{are} linked
(essentially equal).

\begin{dfn}\label{d:link3}
Suppose that $(\lam_1,\fup_1)$ and $(\lam_2, \fup_2)$ are
geolaminations with full portraits. They are said to be
\emph{linked} (\emph{essentially equal}) if there are geolaminations
$(\lamm_1, \qcp_1)$ and $(\lamm_2, \qcp_2)$ with qc-portraits
$\qcp_1\prec \fup_1$, $\qcp_2\prec \fup_2$ such that
$\lamm_1\supset \lam_1$, $\lamm_2\supset \lam_2$, where $(\lamm_1,
\qcp_1)$ and $(\lamm_2, \qcp_2)$ are linked (essentially equal).
\end{dfn}

Let us prove a simple but useful lemma.

\begin{lem}\label{l:actri}
If $\lam$ is a geolamination with an all-critical triangle $T$ then,
for any two full portraits $\fup^1, \fup^2$ of $\lam$, the
geolaminations with full portraits $(\lam, \fup^1)$ and $(\lam,
\fup^2)$ are essentially equal.
\end{lem}

\begin{proof} It suffices to declare $T$ a special critical cluster.
\end{proof}

In Definition~\ref{d:link3} (which is rather general), we require the
existence of geolaminations $\lamm_1$, $\lamm_2$ with qc-portraits
$\qcp_1$, $\qcp_2$ satisfying specific properties. It turns out that
there are convenient (though not as general) ways to comply with
this definition.

\begin{lem}\label{l:link3}
Suppose that $(\lam_1,\fup_1)$ and $(\lam_2, \fup_2)$ are
geolaminations with full portraits. Suppose that $\qcp_1\prec
\fup_1$ and $\qcp_2\prec \fup_2$ are linked $($essentially equal$)$
qc-portraits such that every collapsing quadrilateral of $\qcp_i$,
where  $i=1,2$, shares a pair of opposite edges with the
corresponding set of $\fup_i$, in which it lies.
Then $(\lam_1,\fup_1)$ and $(\lam_2, \fup_2)$ are linked
$($essentially equal$)$.
\end{lem}

\begin{proof}
Let $Q$ be a set of one of our qc-portraits, $Y$ be the
corresponding set of the corresponding full portrait, and $Q\subset
Y$. Then the image of $Q$ is an edge of $\si_3(Y)$ and, therefore, a
leaf of the corresponding geolamination. Hence forward images of
sets of inserted qc-portraits $\qcp_1$, $\qcp_2$ will not generate
linked leaves. Standard arguments show that we can then pull back
sets of $\qcp_1, \qcp_2$ and thus construct the geolaminations
$\lamm_1$ and $\lamm_2$ as required in Definition~\ref{d:link3}.
\end{proof}

The following lemma is a key combinatorial fact about the tags.

\begin{lem}\label{l:interlink}
Suppose that $(\lam_a,\fup_a)$ and $(\lam_x, \fup_x)$ are
geolaminations with full portraits, and $\lam_a$ is dendritic.
If their mixed tags are non-disjoint, then these geolaminations are linked
or essentially equal.
\end{lem}

The proof of Lemma \ref{l:interlink} is mostly non-dynamical and involves
checking a variety of cases.
We split the proof into several propositions.
While a general argument exists when
both geolaminations are bicritical, arguments for one or two
unicritical geolaminations are more complicated due to the number of
different full portraits that can be associated with them.

\begin{prop}\label{p:critrilink}
Suppose that $(\lam_i,\fup_i)$ and $(\lam_j, \fup_j)$ are
geolaminations with full portraits and non-disjoint mixed tags. If
$\lam_i$ contains an all-critical triangle $T$, then either $\lam_i$
and $\lam_j$ are linked $($essentially equal$)$, or $\lam_j$
contains no critical sets of degree three, and two edges of $T$ are
contained in the distinct sets of $\fup_j$.
\end{prop}

\begin{proof}
Let $T$ have vertices $a, b$ and $c$, set $\fup_i=(D_1, D_2)$ and
$\fup_j=(E_1, E_2)$. Then $D_2$ is either $T$ or an edge of $T$; so,
$\si_3(D_2)=\si_3(T)=x\in\si_3(E_2)$. If $\lam_j$ has $T$ as a gap,
then, by Lemma~\ref{l:actri}, the geolaminations with full portraits
$(\lam_i,\fup_i)$ and $(\lam_j,\fup_j)$ are essentially equal. If
$\lam_j$ has a critical gap $G$ of degree three that is not a
critical triangle, then, by Definition~\ref{d:admicri}, we have
$E_1=E_2=G$, and $x\in\si_3(G)$ implies that $T\subset G$. Clearly,
in this case, we can choose equal qc-portraits in $\fup_i$, $\fup_j$,
which again shows that $(\lam_i,\fup_i)$ and $(\lam_j, \fup_j)$ are
essentially equal.

Assume now that $\lam_j$ has no critical set of degree three, so that
$E_1\ne E_2$ are of degree two. Since $x\in \si_3(E_2)$, the set $E_2$
contains, say, $\ol{ab}$ denoted so that $(a, b)$ is a hole
of $T$. Then it is easy to see that the vertices of $E_1$ belong to
$[b,a]$, while the vertices of $\coc(E_1)$ belong to $[a, b]$. Since
$\coc(E_1)$ is non-disjoint from $\coc(D_1)\subset T$ then
$\coc(E_1)$ must contain either $a$ or $b$; let $a\in \coc(E_1)$.
This implies that $\ol{bc}\subset E_1$.
\end{proof}

It is easy to give examples, when the second case from
Proposition~\ref{p:critrilink} is realized. Indeed, let $\lam_1$ be
a geolamination with an all-critical triangle $T$ with vertices $0,
\frac13, \frac23$ such that $\fup_1=(\ol{0\frac13},
\ol{\frac13\frac23})$. Let $\lam_2$ be a geolamination with a leaf
$\ol{0\frac12}$ and critical leaves $\ol{0\frac13}, \ol{0\frac23}$
so that $\fup_2=(\ol{0\frac13}, \ol{0\frac23})$. Either
geolamination can be constructed by means of iterative pull backs of
the already given sets (observe that $T$ maps to $0$ and
$\si_3(\ol{0\frac12})=\ol{0\frac12}$). However, $(\lam_1, \fup_1)$
and $(\lam_2, \fup_2)$ are neither linked nor essentially equal,
because the only qc-portrait $\qcp_2\prec \fup_2$ is
$\qcp_2=\fup_2=(\ol{0\frac13}, \ol{0\frac23})$, any qc-portrait
$\qcp_1\prec \fup_1$ must contain $\ol{\frac13 \frac23}$ as its
second set, and then $\qcp_1$ and $\qcp_2$ cannot be linked or
essentially equal.

However if we make an extra assumption on one of the geolaminations
being dendritic, the second case of Proposition~\ref{p:critrilink}
does not realize.

\begin{lem}\label{l:critrilink}
Suppose that $(\lam_1,\fup_1)$ and $(\lam_2, \fup_2)$ are
geolaminations with full portraits and non-disjoint mixed tags. If
$\lam_1$ is dendritic, and one of the geolaminations $\lam_1$, $\lam_2$
contains an all-critical triangle, then $(\lam_1, \fup_1)$ and
$(\lam_2, \fup_2)$ are linked or essentially equal.
\end{lem}

\begin{proof}
We may assume that the second case of Proposition~\ref{p:critrilink}
holds, and one of geolaminations $\lam_1$, $\lam_2$ (say, $\lam_i$) has
an all-critical triangle $T$, and the other one (say, $\lam_j$) is such that the
full portrait $\fup_j$ has distinct critical sets containing two
distinct edges of $T$. Then the two critical sets of $\fup_j$ are
distinct and non-disjoint. Hence, by Lemmas \ref{l:perfectd} and
\ref{l:cridisj}, the geolamination $\lam_j=\lam_2$ is not dendritic
while $\lam_i=\lam_1$ is dendritic. Replace $\fup_1$ by a full
portrait $\fup'_1$ consisting of the two critical edges of $T$
contained in the distinct critical sets of $\fup_2$. Then by
definition $(\lam_1, \fup'_1)$ is essentially equal to $(\lam_2,
\fup_2)$. Since $\lam_1$ is perfect, this implies that $\lam_1\subset
\lam_2$, and hence $T$ is a gap of $\lam_2$. By Lemma~\ref{l:actri},
we obtain the desired.
\end{proof}

A nice paper by Dierk Schleicher \cite{sch04} contains a full
treatment of the case, when both geolaminations are unicritical and
of degree $d$. We, however, only need a simple fact concerning
unicritical \emph{cubic} geolaminations.

\begin{lem}\label{l:uniclink}
Suppose that $(\lam_1,\fup_1)$ and $(\lam_2, \fup_2)$ are
unicritical geolaminations with full portraits. Assume that mixed
tags of these geolaminations are non-disjoint. Then
$(\lam_1,\fup_1)$ and $(\lam_2, \fup_2)$ are linked or essentially
equal.
\end{lem}

\begin{proof}
If both geolaminations have all-critical triangles, then the claim
follows from Lemma~\ref{l:actri}. If exactly one of the two geolaminations
has an all-critical triangle, then the claim follows from
Proposition \ref{p:critrilink}. Suppose that neither geolamination has an all-critical
triangle, and let their critical sets be $C$ (for $\lam_1$) and
$K$ (for $\lam_2$). If $\si_3(C)\cap \si_3(K)$ contains a point
$x\in \uc$ then the entire all-critical triangle $\si_3^{-1}(x)$ is
contained in $C\cap K$, which by definition implies the claim.
Otherwise, we may assume that an edge $\ovc$ of $\si_3(C)$ crosses
an edge $\ovk$ of $\si_3(K)$. This implies
that the hexagon $\si_3^{-1}(\ovk)=\widetilde K\subset K$ and the
hexagon $\si_3^{-1}(\ovc)=\widetilde C\subset C$ have alternating vertices.
This immediately implies the claim of the lemma too.
Thus, in all possible cases, $(\lam_1,\fup_1)$ and $(\lam_2,
\fup_2)$ are linked or essentially equal.
\end{proof}

We are now ready to prove Lemma \ref{l:interlink}.

\begin{proof}[Proof of Lemma \ref{l:interlink}]
By the preceding results, we may assume that neither geolamination has
a critical triangle and that at least one geolamination is
bicritical.
Set $\fup_a=(C_1,C_2)$ and $\fup_x=(K_1,K_2)$.
Since $\coc(C_1)\cap\coc(K_1)\ne\0$, we may suppose that
either there is a point $z\in \coc(C_1)\cap \coc(K_1)\cap \uc$, or
there are leaves $\ovc=\ol{ab}$, $\ovk=\ol{xy}$ of $\coc(C_1)$,
$\coc(K_1)$ that cross, in which case we assume that $(a,b)$, $(x,y)$
are holes of $\coc(C_1)$, $\coc(K_1)$ and $a\le x\le b\le y$ (the
argument in the case when $x\le a\le y\le b$ is similar).

Consider first the case, when there is a point $z\in \coc(C_1)\cap
\coc(K_1)\cap \uc$. Then $\ol{(z+\frac13)(z+\frac23)}=\n$ is a
critical chord shared by $C_1$ and $K_1$. Now, if $\lam_a$ is
bicritical, then $C_2\cap \uc\subset (z+\frac23, z+\frac13)$ whereas,
if $\lam_a$ is unicritical, then the appropriate part of $C_2=C_1$,
on which $\si_3$ is two-to-one, has vertices belonging to
$[z+\frac23, z+\frac13]$. Similarly, $K_2$ (or the appropriate
two-to-one part of $K_2$ if $\lam_x$ is unicritical) is contained in
$[z+\frac23, z+\frac13]$.

Consider now the sets $\si_3(C_2)$ and
$\si_3(K_2)$. If $\si_3(z)\in \si_3(C_2)\cap \si_3(K_2)$, then at least one of the
points $z+\frac13$, $z+\frac23$ belongs to $C_2\cap C_1$. This
implies, by Lemmas \ref{l:perfectd} and \ref{l:cridisj}, that
$C_1=C_2=C$, and $\lam_a$ is unicritical. Hence, by the assumptions from the beginning of the proof,
$\lam_x$ is bicritical. We claim that there exists
a critical chord with the image $\si_3(z)$, contained in $K_2$ and
not equal to $\n$.
Suppose otherwise. Then $\n$ is shared by $K_1$ and $K_2$, which
implies that one of these sets equals $\n$ while the other one is a
gap $G$ of degree two with more than three vertices located between
$z$ and $\n$ and such that the map $\si_3$ is exactly two-to-one on
it except for the points $z, z+\frac13$ and $z+\frac23$. It follows
from the definition that the co-critical set of this gap cannot
contain $\coc(K_1)$, which implies that in fact $K_2=G$, and the
desired chord exists. Clearly, this chord is shared by $C$ and $K_2$
which implies that $\lam_a$ and $\lam_x$ are essentially equal.

If now $\si_3(C_2)\cap \si_3(K_2)\cap \uc$ contains a point distinct from
$\si_3(z)$, then its pullback to $(z+\frac23, z+\frac13)$ is a
critical chord shared by $C_2$ and $K_2$, again showing that $\lam_a$
and $\lam_x$ are linked or essentially equal. Finally, if $\si_3(C_2)$ and $\si_3(K_2)$
have linked edges then their pullbacks to $[z+\frac23, z+\frac13]$
are strongly linked collapsing \ql s sharing edges with $C_2$ and
$K_2$, respectively. As above, this implies that $\lam_a$ and
$\lam_x$ are linked.

Assume now that there are crossing leaves $\ovc=\ol{ab}$, $\ovk=\ol{xy}$ of
$\coc(C_1)$, $\coc(K_1)$; let $(a,b)$, $(x,y)$ be the holes of
$\coc(C_1)$, $\coc(K_1)$, and 
$a<x<b<y$. We claim that $y\le
a+\frac13$. Indeed, otherwise $[b, a+\frac13]\subset [x, y)$ which
implies that $[\si_3(b), \si_3(a)]\subset [\si_3(x), \si_3(y))$. On
the other hand, by Proposition~\ref{p:colocation}, we have
$\si_3(C_2)\subset [\si_3(b), \si_3(a)]$ and $\si_3(K_2)\subset
[\si_3(y), \si_3(x)]$. Since $\si_3(C_2)\cap \si_3(K_2)\ne \0$, then
in fact $b=x$, a contradiction. Thus, the endpoints of $\ovc$ and
$\ovk$ belong to an interval of length at most $\frac13$. By
Proposition~\ref{p:linkco}, the gaps $C_1$ and $K_1$ contain
strongly linked \ql s $\coc(\ovc)=Q^c_1$ and $\coc(\ovk)=Q^k_1$.

Consider the closure $R^c$ of the component of $C_2\sm Q^c_1$ that
does not map one-to-one onto its image. Considering unicritical and
bicritical cases separately, we see that $\si_3$ maps $R^c$ onto its
image two-to-one. Similarly, we can define the set $R^k$. It follows
that the arc $[x+\frac23, x+\frac13]$ contains vertices of both
$R^c$ and $R^k$ and that $\si_3(C_2)=\si_3(R^c)$ and
$\si_3(K_2)=\si_3(R^k)$. As above, it follows that $R^c$ and $R^k$
contain strongly linked critical generalized \ql s or
critical generalized \ql s which share a spike; these \ql s
share two opposite
sides with $C_2$ and $K_2$, respectively. Hence in this case
$\lam_a$ and $\lam_x$ are linked or essentially equal as well.
\end{proof}

We are ready to prove Theorem~\ref{t:dendrilink}.

\begin{thm}\label{t:dendrilink}
If $(\lam_a,\fup_a)$ and $(\lam_x, \fup_x)$ are cubic
geolaminations with full portraits $\fup_a=(C_1, C_2)$,
$\fup_b=(K_1, K_2)$, and $\lam_a$ is dendritic, then they have
non-disjoint mixed tags if and only if one of the following holds:

\begin{enumerate}

\item $\lam_a$ has an all-critical triangle $T$ as a gap, $\lam_a=\lam_x$,
and it is not true that $C_1$ and $K_1$ are distinct edges of\, $T$;

\item $\lam_a$ does not have an all-critical triangle as a gap,
$\lam_a\subset\lam_x$, and $\fup_a\succ \fup_x$.

\end{enumerate}

\end{thm}

\begin{proof}
Suppose that mixed tags of $(\lam_a, \fup_a)$ and $(\lam_x, \fup_x)$
are non-disjoint. Then, by Lemma~\ref{l:interlink}, the
geolaminations with full portraits $(\lam_a, \fup_a)$ and $(\lam_x,
\fup_x)$ are linked or essentially equal. Since $\lam_a$ is perfect,
by Theorem~\ref{t:noesli}, this implies that always $\lam_x\supset
\lam_x^c=\lam_a$ (recall that $\lam_x^c$ is the maximal perfect
sublamination of $\lam_x$).

Assume that $\lam_a$ has an all-critical triangle $T$ as a gap.
Then, by \cite{bl02, kiw02}, the geolamination $\lam_a$ does not
have wandering polygons and, by \cite{kiw02}, if $G$ is a periodic
gap of $\lam_a$ with more than three vertices, then all its vertices
belong to the same periodic orbit. Hence in this case $\lam_x\supset
\lam_a$ implies that $\lam_x=\lam_a$. Moreover, since mixed tags of
our geolaminations are non-disjoint, $\coc(C_1)\cap \coc(K_1)\ne
\0$. Clearly, of $C_1$ and $K_1$ are distinct edges of $T$ then
$\coc(C_1)\cap \coc(K_1)=\0$. On the other hand, it is easy to
verify (considering a few cases) that otherwise the mixed tags are
non-disjoint as desired. This completes the proof of (1).

On the other hand, assume that $\lam_a$ does not have an
all-critical triangle as a gap. Then since the mixed tags of
$\lam_a$ and $\lam_x$ are non-disjoint and $\lam_a\supset \lam_x$,
then $\fup_a\succ \fup_x$ as desired.

The opposite direction of the theorem follows from definitions.
\end{proof}

Observe that the condition from Theorem~\ref{t:dendrilink}(1) that
it is not true that $C_1$ and $K_1$ are distinct edges of $T$ is
equivalent to the condition that either $C_1\supset K_1$, or
$K_1\supset C_1$.

\subsection{Upper semi-continuous tags}
We will now introduce a topology in the space of tags.

\begin{dfn}\label{d:usc} A collection $\mD=\{D_\alpha\}$ of compact
and disjoint subsets of a metric space $X$ is \emph{upper
semicontinuous $($USC$)$} if, for every $D_\alpha$ and every open set
$U\supset D_\alpha$, there exists an open set $V$ containing
$D_\alpha$ so that for each $D_\beta\in \mD$, if $D_\beta\cap
V\ne\0$, then $D_\beta\subset U$.
\end{dfn}

\begin{thm}[\cite{dave86}]\label{t:dav} If $\mD$ is an
upper se\-mi\-con\-ti\-nuous decomposition of a separable metric
space $X$, then the quotient space $X/\mD$ is also a separable
metric space.
\end{thm}

To apply Theorem~\ref{t:dav}, we need Theorem~\ref{t:tagusc}. However,
first we study limits of finite critical sets of \gl s.

\begin{lem}\label{l:decrilim}
Let $C_1,C_2,\dots$ be a sequence of finite critical sets of \gl s
$\lam_i$ converging to a set $C$. Then $C$ is a critical set $($in
particular, $C$ is not a gap of degree one$)$, and $C$ is not periodic.
\end{lem}

\begin{proof}
We may assume that all sets $C_i$ have degree $k$. Then the degree
of $C$ is at least $k$, and hence $C$ is critical. If $C$ is
periodic of period, say, $n$, then, since it is critical, it is an
infinite gap. Then the fact that $\si_d^n(C)=C$ implies that any gap
$C_i$ sufficiently close to $C$ will have its $\si_d^n$-image also
close to $C$, and therefore coinciding with $C_i$. Thus, $C_i$ is
$\si_d$-periodic, which is impossible because $C_i$ is finite and
critical.
\end{proof}

We will also need the following elementary observation.

\begin{lem}\label{l:crifar}
Suppose that $C_1$, $C_2$ are distinct critical sets of a cubic
geolamination. Then there exists a point in one of them, whose
distance to the other critical set (measured along the circle) is at
least $\frac1{12}$.
\end{lem}

\begin{proof}
Choose a chord $\ell$ separating $C_1\sm C_2$ from $C_2\sm C_1$.
Clearly, there exist two semi-open strips $L$ and $R$ located on
either side of $\ell$, each of which is a convex hull of two
circular arcs of length $\frac1{12}$ sharing endpoints of $\ell$
with one circular chord-edge (not equal to $\ell$) removed. If one
of the critical sets is not contained in $L\cup R$ the claim
follows. Hence we may assume that $C_1\subset L, C_2\subset R$ which
implies that both $L$ and $R$ contain a critical chord, a
contradiction.
\end{proof}

We are now ready to show that Theorem~\ref{t:dav} applies to our
tags.

\begin{thm}\label{t:tagusc}
The family $\{\ta(\zc)\}$ of tags of critical patterns of dendritic
geolaminations forms an upper semicontinuous decomposition of
$\{\ta(\zc)\}^+$.
\end{thm}

\begin{proof} If $(\lam_1,\zc_1)$ and
$(\lam_2,\zc_2)$ are two dendritic geolaminations with critical
patterns, and $\ta(\zc_1)$ and $\ta(\zc_2)$ are non-disjoint, then, by
Lemma~\ref{l:interlink}, we have $(\lam_1, \zc_1)=(\lam_2, \zc_2)$.
Hence the family $\{\ta(\zc)\}$ forms a decomposition of the 
union of tags of all dendritic geolaminations $\{\ta(\zc)\}^+$.

Suppose next that $(\lam_i,\zc_i)$ is a sequence of dendritic
geolaminations with critical patterns  $\zc_i=(C_i^1,C^2_i)$ and
tags $\coc(C^1_i)\times \si_3(C^2_i)$. Suppose that 
there is a limit point of the sequence $\coc(C^1_i)\times
\si_3(C^2_i)$ that belongs to the tag of a dendritic geolamination
$\lam_D$ with a critical pattern $\zc_D=(C^1_D,C^2_D)$. By
\cite{bmov13}, we  may assume that the sequence $\lam_i$ converges
to an invariant geolamination $\lam_\infty$. Then, by
Lemma~\ref{l:decrilim}, the critical sets $C_i^1, C_i^2$ converge to
critical sets $C_\infty^1, C_\infty^2$ of $\lam_\infty$. By
Lemma~\ref{l:crifar}, if $C^1_i\ne C^2_i$ for all sufficiently large
$i$, then $C_\infty^1\ne C_\infty^2$, and $\fup_\infty=(C^1_\infty,
C^2_\infty)$ is a full portrait of $\lam_\infty$. By the assumption,
$\ta(\zc_D)\cap \ta(\fup_\infty)\ne\0$. By
Theorem~\ref{t:dendrilink}, $\lam_D\subset \lam_\infty$ and
$\fup_\infty\prec \zc_D$. 
Hence $\ta(\lam_\infty, \fup_\infty)\subset \ta(\lam_D,
\zc_D)$.
\end{proof}

Denote the quotient space $\{\ta(\zc)\}^+/\{\ta(\zc)\}$ by $\cml$.
We show that $\cml$ can be viewed as a combinatorial model for
$\cmd_3$. By Theorem~\ref{t:dav}, the topological space $\cml$ is separable and metric. We
denote the quotient map from $\{\ta(\zc)\}^+$ to $\cml$ by
$\pi_{\ta}$. By Corollary~\ref{c:crista}, the map $\Psi_3$ maps a
critically marked dendritic polynomial $(P, \cm)$ to a dendritic
geolamination with a critical pattern $(\lam_P, \zc)$. Together with
our definitions, this implies the following theorem.

\begin{thm}\label{t:polytags}
The composition
$$\Phi_\ta(P, CM)=\pi_\ta\circ \ta \circ \Psi_3(P,
\cm)$$ is a continuous surjective map $\Phi_\ta: \cmd_3\to \cml$.
\end{thm}

\begin{proof}
Let $(P_i, \cm_i)\to (P, \cm)$ with $(P_i, \cm_i)\in \cmd_3$, $(P,
\cm)\in \cmd_3$ and $\Psi_3(P, \cm_i)=(\lam_{P_i}, \zc_i),$ $\Psi_3(P,
\cm)=(\lam_P, \zc)$ with critical sets $C_1,C_2$. Without loss of
generality, we may assume that $(\lam_{P_i}, \zc_i)$ converge in the
Hausdorff sense to $(\lam^\infty, (C_1^\infty, C^\infty_2))$. Then, by
Corollary~\ref{c:crista}, we have $\lam^\infty\supset \lam_P$ and
$C^\infty_i\subset C_i$ for $i=1,2$. By definition, this implies the
desired.
\end{proof}

\end{document}